# ASYMPTOTIC EFFICIENCY AND FINITE-SAMPLE PROPERTIES OF THE GENERALIZED PROFILING ESTIMATION OF PARAMETERS IN ORDINARY DIFFERENTIAL EQUATIONS[1]

By Xin Qi and Hongyu Zhao

*Yale University*

Ordinary differential equations (ODEs) are commonly used to model dynamic behavior of a system. Because many parameters are unknown and have to be estimated from the observed data, there is growing interest in statistics to develop efficient estimation procedures for these parameters. Among the proposed methods in the literature, the generalized profiling estimation method developed by Ramsay and colleagues is particularly promising for its computational efficiency and good performance. In this approach, the ODE solution is approximated with a linear combination of basis functions. The coefficients of the basis functions are estimated by a penalized smoothing procedure with an ODE-defined penalty. However, the statistical properties of this procedure are not known. In this paper, we first give an upper bound on the uniform norm of the difference between the true solutions and their approximations. Then we use this bound to prove the consistency and asymptotic normality of this estimation procedure. We show that the asymptotic covariance matrix is the same as that of the maximum likelihood estimation. Therefore, this procedure is asymptotically efficient. For a fixed sample and fixed basis functions, we study the limiting behavior of the approximation when the smoothing parameter tends to infinity. We propose an algorithm to choose the smoothing parameters and a method to compute the deviation of the spline approximation from solution without solving the ODEs.

**1. Introduction.** Ordinary differential equations (ODEs) are often used to model dynamic processes in engineering, biology and many other areas.

---

Received March 2009; revised June 2009.

[1]Supported in part by NIH Grants P30 DA018343, R01 GM59507 and NSF Grant DMS-07-14817.

For example, the dynamic behavior of gene regulation networks can be modeled by a set of ODEs (see Gardner et al. [9] and Cao and Zhao [5]). These ODEs usually involve many unknown parameters. Ideally, we hope that we can estimate these unknown parameters by some classical parametric estimators, such as least squares estimators or maximum likelihood estimators (MLE). However, most nonlinear ODE systems do not have analytical solutions whereas numerical solutions can be time consuming and it is nontrivial to estimate their values from the observed data that are often very noisy.

Because of the importance of this problem, many methods have been proposed to estimate parameters in ODEs that cannot be solved analytically. One method is through nonlinear least squares (NLS). In this approach, a numerical method is used to approximate the solution of the ODEs at a given trial set of parameter values and initial conditions. The fitted values are input into the nonlinear least squares procedure to update parameter estimates. This NLS approach is computationally intensive since a numerical approximation to the solutions is required for each update of the parameters and initial conditions. In addition, the inaccuracy of the numerical approximation can be a problem, especially for stiff systems.

Another approach, called collocation methods, approximates the solution by a basis function expansion, such as the cubic spline function. Varah [27] suggested a two-stage procedure where the first step fits the observed data by least squares using cubic spline functions without considering the ODEs, and the second step estimates the parameters by least squares solution of the differential equations sampled at a set of points. This approach works well for the simple equations considered, but considerable care is required in the smoothing step and all the variables in the system need to be measured. Ramsay and Silverman [22] and Poyton et al. [20] further developed Varah's method by proposing an iterated principal differential analysis, which converged quickly to the estimates of both the solution and the parameters and had substantially improved bias and precision. However, their approach is a joint estimation procedure in the sense that it optimizes a single roughness-penalized fitting criterion with respect to both the coefficients of the basis expansion and the parameters. The effect of the nuisance parameters on the fit of the model cannot be controlled. For other collocation methods, see Tjoa and Biegler [24], Arora and Biegler [2] and Bock [3].

Most recently, Ramsay et al. [21] proposed a new collocation method called generalized profiling procedure. In this approach, the ODE solution is approximated by a linear combination of basis functions. However, the coefficients of the basis functions are estimated by a penalized smoothing procedure with an ODE-defined penalty. The smoothing parameter controls the trade off between fitting the data with the basis functions and fidelity of the basis functions to the ODEs. Their method has several unique aspects. The computation load is much lower than other methods because it avoids numerically solving ODEs. It can estimate some ODE components even if they are not observed. It is easy to estimate uncertainties in parameter estimates and simulation experiments suggested that there is good agreement between estimated uncertainties and actual estimation accuracies. In addition, this approach does not require a formulation of the dynamic model as an initial value problem in situations where initial values are not available.

Despite these attractive features, little is known about the statistical properties of the estimates from this procedure, such as consistency and asymptotic normality. Furthermore, it is not clear how to choose the smoothing parameters automatically. In this article, to study the asymptotic properties we firstly derive an upper bound on the uniform norm of the difference between the ODE solutions and their approximations in terms of the smoothing parameters and the distance between the approximation space and the solutions (the distance can be controlled by the knots when the cubic spline functions are used as approximations). Then this bound is used to prove the asymptotic consistency of the parameter estimation if the smoothing parameter goes to infinity and the distance between the approximation space and the space of the ODE solutions goes to zero. If the smoothing parameter and the distance satisfy certain conditions on the convergence rate, we prove the asymptotic normality for the parameter estimation and show that its asymptotic covariance matrix is the same as that of the maximum likelihood estimation. Therefore, the profiling procedure is asymptotically efficient. We note that our asymptotic results are also true for partially observed systems (only parts of the components are observed). According to these results, we propose an algorithm to choose the smoothing parameters automatically.

One innovative feature of the profiling procedure is that it incorporates a penalty term to estimate the coefficients in the first step. This penalty is the $L^2$-norm of the difference between the two sides of ODEs which are evaluated by plugging in the approximation functions. From the theory of differential equations, for such penalty (even the $L^\infty$-norm), the bound on the difference between the approximations and the solutions will grow exponentially however small the penalty is. As a result, if the time interval is large, the bound will be too large to be useful. However, the results in Ramsay et al. [21] and our simulation studies indicate that when the smoothing parameter becomes large, the approximations to the solutions are very good. There is no trend of exponentially growing. To explain this phenomena, we fix the sample and the approximation space, and study the limiting situation as the smoothing parameter goes to infinity. We show that any such sequence will have a subsequence converging to one of the minimum functions of the penalty in the approximation space. We study the properties of these minimum functions and give a bound on the uniform norms of the differences between these functions and the solutions in the one-dimensional case and the B-spline bases. The bound depends almost linearly on the length of the time interval or almost does not depend on length of the time interval if we put stronger conditions. This result explains the above noted phenomena and motivates us a method to compute the deviation of the spline approximation from solutions without solving the ODEs.

Olhede [19] outlined some asymptotic results for the profiling procedure proposed by Ramsay and colleagues. In order to achieve asymptotic consistency, Olhede [19] took the number of the B-spline basis functions to be of order $O(n)$ where $n$ is the sample size. Then she imposted the conditions that the penalties have order $O(n^{-\delta})$ for some $\delta > 0$ and the sum of the scaled penalties by smoothing parameters has order $O(n)$. It was derived that the smoothing parameters have order $O(n^{1+\delta})$. However, it is not easy to tune these parameters to ensure that the penalties satisfy the imposed conditions, because only the number of the bases and the smoothing parameters can be tuned and the values of the penalties depend on these two sets of tuning parameters through complex relationships. Therefore, it is better to impose conditions only on the tuning parameters. Furthermore, although we obtain the approximations to the solutions by choosing a set of basis functions and computing the coefficients by solving a penalized optimization problem, the approximations only depend on the linear space spanned by the bases. In fact, if we choose another set of bases in the same space, we should get the same approximations. Hence, an appropriate theory should put conditions on the approximation space instead of the bases. In her discussion, Olhede [19] proposed to use $L^\infty$-penalty instead of $L^2$ or $L^1$-penalty. Our asymptotic results hold for all these penalties, although the smoothing parameters take different convergence rates for different norms. Unfortunately, the bound that is exponentially increasing in time in Theorem 3.1 will stay for all the penalties. The results for fixed sample where the smoothing parameters converge to infinity can only be shown for $L^2$-penalty because there we use the special property of $L^2$-norm that we can change the integrals for time and the differentiation with respect to the parameters under this norm.

Lele [16] raised the question as to what kind of asymptotics is appropriate: infill asymptotics, or increasing domain asymptotics, or both. Here we choose the infill asymptotics for the following reasons. First, one key point in our proof is the uniform boundedness of the solutions to ODEs on a finite time interval for all the parameters in a compact subset of the parameter space. In order to do that, we make assumptions about the existence of the solutions and the smoothness of the functions in ODEs. However, the existence does not always hold or is not easy to check for nonlinear ODEs [see Remark 2(1) after Assumption 2]. If we choose a fixed time interval, for example, the largest sampled time as the endpoint of the interval, there exists at least a neighborhood of the true parameters such that the solutions exist for the parameters in this neighborhood. However, these conditions and assumptions could be seriously violated for increasing domain asymptotics. Second, the difference between the basis function approximations and the true solutions could increase exponentially with time (see Theorem 3.1 and the remarks after it). This will affect the accuracy of the estimate. Last, if we are only interested in estimating the parameters, it is adequate to sample an increasing number of points in a finite interval under the assumptions that the model is correct, the functions in ODEs are smooth, and the parameters are identifiable.

Now, we describe the model and the profiling procedure in detail. Let the parameter space $\Theta$ be an open and convex subset of $\mathbb{R}^d$. We use $\theta_0$ to denote the true parameter. Consider the following ODEs:

$$
\begin{aligned}
\frac{d\mathbf{x}}{dt}(t) &= F(\mathbf{x}(t), \mathbf{z}(t), t, \theta), \\
\frac{d\mathbf{z}}{dt}(t) &= G(\mathbf{x}(t), \mathbf{z}(t), t, \theta), \\
\mathbf{x}(0) &= x_0, \qquad \mathbf{z}(0) = z_0,
\end{aligned}
\tag{1.1}
$$

on time interval $[0, T]$ with $\mathbf{x} : [0, T] \to \mathbb{R}^{d_1}$, $\mathbf{z} : [0, T] \to \mathbb{R}^{d_2}$, $F : \mathbb{R}^{d_1} \times \mathbb{R}^{d_2} \times [0, T] \times \Theta \to \mathbb{R}^{d_1}$, and $G : \mathbb{R}^{d_1} \times \mathbb{R}^{d_2} \times [0, T] \times \Theta \to \mathbb{R}^{d_2}$. $F$ and $G$ have known functional forms with some unknown parameters, and $x_0$ and $z_0$ are the initial values in $\mathbb{R}^{d_1}$ and $\mathbb{R}^{d_2}$. Suppose that the initial values can be chosen from a convex open region $\Gamma \in \mathbb{R}^{d_1} \times \mathbb{R}^{d_2}$. We assume that for each $\theta$, the initial value problem (1.1) has a unique solution $(\mathbf{x}(\theta, t), \mathbf{z}(\theta, t))$ on $[0, T]$. We use the bold face letters to denote the functions on $[0, T]$.

The following is a concrete example from Ramsay et al. [21]. Consider the FitzHugh–Nagumo equations which describe the behavior of spike potentials in the giant axon of squid neurons:

$$
\begin{aligned}
\frac{d\mathbf{V}}{dt}(t) &= c\left(\mathbf{V} - \frac{\mathbf{V}^3}{3} + \mathbf{R}\right), \\
\frac{d\mathbf{R}}{dt}(t) &= -\frac{1}{c}(\mathbf{V} - a + b\mathbf{R}),
\end{aligned}
$$

where $\mathbf{V}$ is the voltage across an axon membrane and $\mathbf{R}$ is a recovery variable summarizing outward currents. The parameters are $\theta = (a, b, c)$, and the time interval is $[0, 20]$.

Assuming the underlying model (1.1), suppose that $(Y_1, T_1), \ldots, (Y_n, T_n)$ are i.i.d. observations, where the $T_i$ is the sample time and $Y_i$ is the observed data at time $T_i$. $T_i$'s are independent random variables distributed on interval $[0, T]$ with distribution $Q$ and $Y_i$'s take the values in $\mathbb{R}^{d_3}$. To fit a model of form (1.1) to the observed data suppose we have the following data fitting criterion

$$
H_n(\mathbf{x}(\theta, \cdot)) = -\frac{1}{n} \sum_{i=1}^{n} g(Y_i, \mathbf{x}(\theta, T_i)),
$$

where $g(y, x)$ is a function on $\mathbb{R}^{d_3} \times \mathbb{R}^{d_1}$ and $H_n$ is a functional on the space of the function on $[0, T]$ such that for any function $\mathbf{x}(t)$,

$$H_n(\mathbf{x}) = -\frac{1}{n} \sum_{i=1}^{n} g(Y_i, \mathbf{x}(T_i)).$$

Here, we assume that only part of the systems are observable, e.g., $Y_i$'s only depend on the first $d_1$ components of the solution. The following are two such examples.

EXAMPLE 1 (Nonlinear least squares).   Suppose that $d_1 = 1$ and

$$Y_i = \mathbf{x}(\theta_0, T_i) + \varepsilon_i,$$

where $\{\varepsilon_i, i = 1, \ldots, n\}$ are independent random variables with the same distributions $N(0, \sigma^2)$. Here we take $g(y, x) = (y - x)^2$.

EXAMPLE 2 (Logistic regression).   Suppose that $d_1 = 1$ and the conditional distribution $Y_i | T_i = t$ is Bernoulli with success probability

$$p(\theta_0, t) = \frac{e^{\mathbf{x}(\theta_0, t)}}{1 + e^{\mathbf{x}(\theta_0, t)}}.$$

Hence,

$$g(y, x) = \log(1 + e^x) - xy.$$

For simplicity, we shall restrict ourselves to the case $d_1 = d_2 = d_3 = 1$. First, we introduce some function spaces and some norms in those spaces. Consider the space of continuously differentiable functions

$$C^1([0, T]) = \{f : \text{both } f \text{ and } f' \text{ are continuous functions on } [0, T]\}$$

and the space of square integrable functions

$$L^2([0, T]) = \left\{ f : f \text{ is a measurable function on } [0, T] \text{ and } \int_0^T |f(t)|^2 \, dt < \infty \right\}.$$

We mainly consider the functions in the space $C^1[0, T]$. For any $f \in C^1[0, T]$, define

$$\|f\|_\infty = \sup_{t \in [0, T]} \{|f(t)|\},$$

$$\|f\|_{L^2([0, T])} = \left[ \int_0^T |f(t)|^2 \, dt \right]^{1/2}.$$

We have two inequalities for these two norms which we will use below.

$$(1.2) \qquad\qquad \|f\|_{L^2([0, T])}^2 \le T \|f\|_\infty^2$$

and by Hölder inequality,

$$\left| \int_0^t f(s)\,ds \right|^2 \leq t \int_0^t |f(s)|^2\,ds \leq T \int_0^T |f(s)|^2\,ds \qquad \forall t \in [0, T],$$

we have

$$(1.3) \qquad \left\| \int_0^t f(s)\,ds \right\|_\infty^2 \leq T \|f(s)\|_{L^2([0,T])}^2,$$

where $\left\| \int_0^t f(s)\,ds \right\|_\infty$ means the norm of the function $h(t) = \int_0^t f(s)\,ds$, $0 \leq t \leq T$.

In this paper, we consider the following penalty, for any $\theta$ and $\mathbf{x}, \mathbf{z} \in C^1[0, T]$,

$$
\begin{aligned}
J(\mathbf{x}, \mathbf{z}, \theta) &= \int_0^T \left[ \frac{d\mathbf{x}}{dt}(t) - F(\mathbf{x}(t), \mathbf{z}(t), t, \theta) \right]^2 dt \\
(1.4) \qquad &\quad + \int_0^T \left[ \frac{d\mathbf{z}}{dt}(t) - G(\mathbf{x}(t), \mathbf{z}(t), t, \theta) \right]^2 dt \\
&= \left\| \frac{d\mathbf{x}}{dt} - F(\mathbf{x}, \mathbf{z}, t, \theta) \right\|_{L^2[0,T]}^2 + \left\| \frac{d\mathbf{z}}{dt} - G(\mathbf{x}, \mathbf{z}, t, \theta) \right\|_{L^2[0,T]}^2.
\end{aligned}
$$

Note that $J$ was used to denote the whole penalized log likelihood criterion in Ramsay et al. [21] and this is different from our definition. We use $\mathbb{P}_n$ to denote the empirical measure. For example,

$$\mathbb{P}_n g(Y, \mathbf{x}(\cdot)) = \frac{1}{n} \sum_{i=1}^n g(Y_i, \mathbf{x}(T_i)).$$

Suppose that $\{\mathbb{L}_n, n \geq 1\}$ is a sequence of finite-dimensional subspaces of $C^1[0, T]$. We will use the functions in $\mathbb{L}_n$ to approximate the solutions of ODEs. For example, we can choose $\mathbb{L}_n$ to be the space of cubic spline functions with knots $\tau^{(n)} = (0 = t_1^{(n)} < \cdots < t_{k_n}^{(n)} = T)$. Define

$$|\tau^{(n)}| = \max_{2 \leq i \leq k_n} |t_i - t_{i-1}|.$$

Let $|\tau^{(n)}| \to 0$ as $n \to \infty$.

Sometimes the initial values of the systems are unknown. In this case, we have to regard the initial values as nuisance parameters. Define $\theta^* = (\theta, x, z)$, the combination of the parameters and the initial values. Let $\theta_0^*$ be the combination of the true parameters and the true initial values. We rewrite the solutions of (1.1) as

$$(\mathbf{x}(\theta^*, t), \mathbf{z}(\theta^*, t)).$$

In the next section, we describe the profiling procedure in details and propose a method of the choice of the smoothing parameter. In Section 3, the consistency and the asymptotic efficiency are given. In Section 4, we study the limit behavior of the basis function approximations to the true solutions to ODEs as the smoothing parameter goes to infinity. The proofs of these results are given in the last section.

**2. The generalized profiling procedure and the choice of the smoothing parameters.** In the generalized profiling procedure proposed in Ramsay et al. [21], a finite-dimensional space $\mathbb{L}_n$ of functions in $[0, T]$ is chosen firstly (this is equivalent to choosing a set of basis functions). Approximations functions to the solutions of the ODEs will be chosen from $\mathbb{L}_n$ (this is equivalent to choosing the coefficients of the basis functions). One innovative feature of this procedure is that the approximations are chosen by solving the following penalized optimization problem: for any $\theta^* = (\theta, x, z) \in \Theta \times \Gamma$,

$$\text{maximize } H_n(\hat{\mathbf{x}}) - \lambda J(\hat{\mathbf{x}}, \hat{\mathbf{z}}, \theta),$$
(2.1)
$$\text{subject to } \hat{\mathbf{x}} \in \mathbb{L}_n, \hat{\mathbf{z}} \in \mathbb{L}_n, \hat{\mathbf{x}}(0) = x, \hat{\mathbf{z}}(0) = z,$$

where the penalty $J$ regularizes the approximations by controlling the size of the extent that the approximation functions fail to satisfy the ODEs and the $\lambda$ is the smoothing parameter which controls the amount of regularization. In this paper, to simplify notation, we choose the same $\lambda$ for all the components in the penalty. Ramsay et al. [21] allowed the different smoothing parameters for different components. The asymptotic results in this paper can be easily extended to this latter case that the smoothing parameters take different values if all the smoothing parameters have the same asymptotic orders. The existence of global solutions to (2.1) will be discussed below in this section. Here, we assume that the global solutions exist. Let $(\hat{\mathbf{x}}(\theta^*, \lambda, t), \hat{\mathbf{z}}(\theta^*, \lambda, t))$ be solutions to (2.1). They depend on both $\lambda$ and $\theta^*$. Then we plug them to the functional $H_n$. The estimates $\hat{\theta}^*(\lambda)$ are obtained by maximizing $H_n(\hat{\mathbf{x}}(\theta^*, \lambda, \cdot))$ with respect to $\theta^*$. Small $\lambda$ makes both the optimization problems of (2.1) and maximizing $H_n(\hat{\mathbf{x}}(\theta^*, \lambda, \cdot))$ robust with respect to the poor initial guesses, but gives bad approximations to the solutions. On the other hand, large $\lambda$ gives good approximations, but the optimization problems are sensitive to the initial values.

In Theorem 3.1 below, the uniform norm of the differences between the exact solutions and the basis function approximations are bounded by a sum of two terms, $O_p(1/\sqrt{\lambda})$ and $O(r_n)$, where $r_n$ is some kind of distance between $\mathbb{L}_n$ and the solutions to ODEs. $r_n$ does not depend on $\lambda$. When $\lambda$ becomes very large, the bound will be dominated by the second term. In this case, it is useless to increase the $\lambda$. We can use these to explain the patterns of Figure 6 in Ramsay et al. [21] and Figure 15 in Huang [12]. Those pictures have different knots and observations. But they have the similar patterns. When we increase the $\lambda$, at first the parameter estimates become better (both of the bias and the variance). But after some point, increasing the $\lambda$ has little effect on the parameter estimates.

The above discussion suggests that we should initially choose a small $\lambda$, then we increase the $\lambda$ until the parameter estimates become stable. In Section 2.8.1 in Ramsay et al. [21], the authors proposed to stop increasing $\lambda$ when the norm of the difference between the solutions to ODEs and the approximations begins increase after obtaining a minimum. Because the penalties are nonlinear functionals of the approximation functions, the approximation functions may depend on $\lambda$ in a complex way. There may exist many local minima before $\lambda$ becomes large enough. Hence, it may be better to design the stopping rule according to the performance of the parameter estimates. One of the major advances of this paper is the fact that it provides accurate estimations of the confidence intervals of the estimated parameters. We can compare the confidence intervals for different $\lambda$ to decide whether we increase $\lambda$ or stop. But computing confidence intervals is time consuming. We can start calculating the confidence intervals when $\lambda$ is moderately large. We propose the following algorithm:

(1) Choose the space $\mathbb{L}_n$, a moderately large positive number $\lambda_0$ and a small number $\alpha$.

(2) Choose a small initial value for $\lambda$, and a initial guess for $\theta^*$.

(3) Obtain the estimates $\hat{\theta}^*(\lambda)$ by maximizing $H_n(\hat{\mathbf{x}}(\theta^*, \lambda, \cdot))$ with respect to $\theta^*$, where for each $\theta^* = (\theta, x, z) \in \Theta \times \Gamma$,

$$(\hat{\mathbf{x}}(\theta^*, \lambda, t), \hat{\mathbf{z}}(\theta^*, \lambda, t)) \in \underset{\substack{\tilde{\mathbf{x}}, \tilde{\mathbf{z}} \in \mathbb{L}_n \\ \tilde{\mathbf{x}}(0) = x, \tilde{\mathbf{z}}(0) = z}}{\arg\max} \ [H_n(\tilde{\mathbf{x}}) - \lambda J(\tilde{\mathbf{x}}, \tilde{\mathbf{z}}, \theta)].$$

(4) If $\lambda < \lambda_0$, set $\hat{\theta}^*(\lambda)$ to be the initial value for next iteration and let $\lambda = 10 \times \lambda$. Go to step (3).

   If $\lambda \geq \lambda_0$, calculate the confidence intervals for the parameter estimations. Compare the confidence intervals with those for previous value of $\lambda$ (if they exist).

   – If the ratios of the overlaps to both of the intervals are larger than $1 - \alpha$, stop and go to step (5).
   – Otherwise, set $\hat{\theta}^*(\lambda)$ to be the initial value for next iteration and let $\lambda = 10 \times \lambda$, then go to step (3).

(5) Set $\lambda_n = \lambda$ and $\hat{\theta}_n^* = \hat{\theta}^*(\lambda)$. $\lambda_n$ is the final choice of the smoothing parameter and $\hat{\theta}_n^*$ is the profiling estimators for the unknown parameters. Set $\hat{\mathbf{x}}_n(\theta^*, t) = \hat{\mathbf{x}}(\theta^*, \lambda_n, t)$ and $\hat{\mathbf{z}}_n(\theta^*, t) = \hat{\mathbf{z}}(\theta^*, \lambda_n, t)$.

Although the penalized optimization problem (2.1) is solved in the finite-dimensional space $\mathbb{L}_n$, the existence of global solutions to (2.1) cannot be easily verified due to the nonlinearity of $J$. However, if we use the norm $\|\cdot\|_\infty$ in $\mathbb{L}_n$, then under our main assumptions, both $J(\hat{\mathbf{x}}, \hat{\mathbf{z}}, \theta)$ and $H_n(\hat{\mathbf{x}})$ are continuous functions of $(\hat{\mathbf{x}}, \hat{\mathbf{z}})$. In this case, the local solutions always exist if we solve (2.1) in any bounded and closed subset of $(\mathbb{L}, \|\cdot\|_\infty)$. In this paper, one important assumption is that the $\hat{\theta}_n^*$ are uniformly tight, that is, given any probability arbitrarily close to 1, one can find a compact subset of the parameter space such that all the $\hat{\theta}_n^*$ belong to that compact subset with this high probability. The solutions to ODEs for any parameter in that compact subset are uniformly bounded by a positive number, say $K$. In this case, if we solve (2.1) in the subset of the functions with norms less than or equal to $K+1$, then the solutions exist and all proofs of the main results are still true by using these local solutions instead of the global solutions. In practice, one can only search a bounded region of the parameter space, so the estimates are uniformly compact. One can solve (2.1) in a large but bounded subset of $\mathbb{L}_n$.

The above algorithm is based on the assumption that the dynamic models are correctly specified. The algorithm in Section 2.8.2 in Ramsay et al. [21] may produce better approximations to the solutions to ODEs for misspecified models. In this case, $\lambda$ may play a slightly different role. Ionides [13] discussed the problem of model misspecification, especially when there are noises in both measurement and dynamics. In this case, it is difficult to define the true trajectories. There are some alternative procedures to deal with this problem, for example, the iterated filtering (see Ionides, Breto and King [14] and Breto et al. [4]) in the frequentist domain or the Bayesian sequential Monte Carlo (see Liu and West [17]) in the Bayesian domain.

**3. Consistency and asymptotic normality.** We now state our main assumptions.

ASSUMPTION 1. $Q$ has a density $f(t)$ with respect to Lebesgue measure on $[0, T]$ and $c \leq f(t) \leq C$ for all $t \in [0, T]$, where $c$ and $C$ are two positive numbers.

REMARK 1. This assumption guarantees that the samples can be taken anywhere and will not be over concentrated on a subset of the time interval.

ASSUMPTION 2. $F, G \in C^3(\mathbb{R} \times \mathbb{R} \times [0, T] \times \Theta)$. For each $\theta \in \Theta$ and each pair of initial values $(x, z) \in \Gamma \subset \mathbb{R} \times \mathbb{R}$, there exists a unique solution $(\mathbf{x}(\theta^*, \cdot), \mathbf{z}(\theta^*, \cdot))$ of the (1.1) on $[0, T]$ and for any $\theta^* = (\theta, x, z) \neq \theta^{*\prime} = (\theta', x', z')$, we have $\mathbf{x}(\theta^*, \cdot) \neq \mathbf{x}(\theta^{*\prime}, \cdot)$

Remark 2.

(1) Here, for developing our theory, we assume that the solutions exist for all the parameters in the parameter space in $[0, T]$. One key point in our paper is that for any compact subset in the parameter space, the solutions for the parameters in this subset are uniformly bounded in $[0, T]$, so the existence of the solutions for all parameters in $[0, T]$ is necessary for our analysis. But it is not always true nor easy to check this assumption for nonlinear ordinary differential equations. Here is a simple example, suppose that $\Theta = (0, \infty)$ and consider the ODE

$$\frac{dx}{dt} = 1 + \theta x^2, \qquad x(0) = 0.$$

For fixed $\theta$, the solution only exists on $[0, \frac{\pi}{2\sqrt{\theta}}]$ and

$$x(t) = \tan \sqrt{\theta} t \qquad \forall t \in \left[0, \frac{\pi}{2\sqrt{\theta}}\right].$$

Hence, for any $T$, if $\theta \geq (\frac{\pi}{2T})^2$, there does not exist any solution in $[0, T]$ with $x(0) = 0$.

In practice, we usually let $T$ be equal to the largest sampled time point. Hence, at least for the true parameter, the solutions exist in $[0, T]$. Since here we assume the smoothness of $F$ and $G$, by Theorem 7.4 in Chapter 1 of Coddington and Levinson [7], there are a neighborhood of the true initial value and a neighborhood of the true parameter such that the solutions exist in $[0, T]$ for the initial values and the parameters in those neighborhoods. If the estimates belong to that neighborhood, all of our results still hold even without the assumption about the existence of the solutions in the whole parameter space. For any parameter for which the solutions do not exist, according to the extension theorem of the solutions to ODEs (see Theorem 3.1 in Hartman [11] or Theorem 1.186 in Chicone [6]), there is a $0 < T' \leq T$ ($T'$ is parameter dependent) such that the solutions exist in $[0, T')$ and the solutions will become unbounded when $t \to T'$. So one can see that when the sample size and the smoothing parameter $\lambda$ are large enough, that parameter cannot be our estimate. But our results may not be true in the situation that our estimates take the values for which the solution do not exists on $[0, T]$ since the solutions are unbounded.

(2) The uniqueness assumption is necessary in our theory. In general, it is not easy to check this assumption for a given ODE and a time interval. However, under the assumption of $F, G \in C^3(\mathbb{R} \times \mathbb{R} \times [0, T] \times \Theta)$, the existence of the solutions of the initial value problems (1.1) in $[0, T]$ is sufficient to guarantee the uniqueness of the solutions in $[0, T]$.

(3) The latter requirement in Assumption 2 means that the parameter estimation problem is identifiable. The model identifiability should be carefully studied before any statistical inference can be made. There is a substantial literature on this issue, for example, Xia [28], Xia and Moog [29], Jeffrey and Xia [15] and Miao et al. [18] investigated the identifiability of HIV dynamic models. There remain many unsolved problems in this area, but they are beyond the scope of this paper.

LEMMA 1. *Under Assumption 2, there exist a sequence of finite-dimensional subspaces $\{\mathbb{L}_n, n \geq 1\}$ of $C^1[0,T]$ such that for any compact subset $\Theta_0$ of $\Theta$ and any compact subset $\Gamma_0$ of $\Gamma$, we have*

$$\lim_{n \to \infty} \sup_{\theta^* \in \Theta_0 \times \Gamma_0} \inf_{\mathbf{w} \in \mathbb{L}_n, \mathbf{w}(0)=x} \left[ \|\mathbf{x}(\theta^*, \cdot) - \mathbf{w}\|_\infty \vee \left\| \frac{d\mathbf{x}}{dt}(\theta^*, \cdot) - \frac{d\mathbf{w}}{dt} \right\|_\infty \right] = 0,$$

$$\lim_{n \to \infty} \sup_{\theta^* \in \Theta_0 \times \Gamma_0} \inf_{\mathbf{v} \in \mathbb{L}_n, \mathbf{v}(0)=z} \left[ \|\mathbf{z}(\theta^*, \cdot) - \mathbf{v}\|_\infty \vee \left\| \frac{d\mathbf{z}}{dt}(\theta^*, \cdot) - \frac{d\mathbf{v}}{dt} \right\|_\infty \right] = 0,$$

*where $\theta^* = (\theta, x, z)$ and $a \vee b$ denotes $max(a, b)$ for any real numbers $a$ and $b$.*

PROOF. Note that in a Euclidean space, a compact subset is just a bounded closed subset. Let $\mathbb{L}_n$ be the space of cubic spline functions with knots $\tau^{(n)} = (0 = t_1^{(n)} < \cdots < t_{k_n}^{(n)} = T)$. Suppose that $|\tau^{(n)}| \to 0$ as $n \to \infty$. Under Assumption 2, $F, G$ have the continuous partial derivatives of third order. Hence, $\frac{d^4\mathbf{x}}{dt^4}(\theta^*, t)$ and $\frac{d^4\mathbf{z}}{dt^4}(\theta^*, t)$ are continuous functions of $(\theta^*, t)$. Because $\Theta_0 \times \Gamma_0$ is a compact set, we have

$$\sup_{\theta^* \in \Theta_0 \times \Gamma_0} \left\| \frac{d^4\mathbf{x}}{dt^4}(\theta^*, \cdot) \right\|_\infty < \infty,$$

$$\sup_{\theta^* \in \Theta_0 \times \Gamma_0} \left\| \frac{d^4\mathbf{z}}{dt^4}(\theta^*, \cdot) \right\|_\infty < \infty.$$

By Theorems 2 and 4 in Hall and Meyer [10],

$$\sup_{\theta^* \in \Theta_0 \times \Gamma_0} \inf_{\mathbf{w} \in \mathbb{L}_n, \mathbf{w}(0)=x} \|\mathbf{x}(\theta^*, \cdot) - \mathbf{w}\|_\infty$$

$$\leq C_0 \sup_{\theta^* \in \Theta_0 \times \Gamma_0} \left\| \frac{d^4\mathbf{x}}{dt^4}(\theta^*, \cdot) \right\|_\infty |\tau^{(n)}|^4 \to 0,$$

$$\sup_{\theta^* \in \Theta_0 \times \Gamma_0} \inf_{\mathbf{w} \in \mathbb{L}_n, \mathbf{w}(0)=x} \left\| \frac{d\mathbf{x}}{dt}(\theta^*, \cdot) - \frac{d\mathbf{w}}{dt} \right\|_\infty$$

$$\leq C_1 \sup_{\theta^* \in \Theta_0 \times \Gamma_0} \left\| \frac{d^4\mathbf{x}}{dt^4}(\theta^*, \cdot) \right\|_\infty |\tau^{(n)}|^3 \to 0,$$

where $C_0 = \frac{5}{384}, C_1 = \frac{9+\sqrt{3}}{216}$. Similarly, we can prove the result for $\mathbf{z}(\theta^*, \cdot)$. $\square$

Let $P_{\theta_0^*}$ be the joint distribution of $\{(Y_1, T_1), \ldots, (Y_n, T_n)\}$, which corresponds to the true parameters and the true initial values. Define function

$$M(\theta^*) = -E_{\theta_0^*}[g(Y_i, \mathbf{x}(\theta^*, T_i))], \qquad \theta^* \in \Theta \times \Gamma,$$

where $E_{\theta_0^*}[\cdot]$ is the expectation with respect to $P_{\theta_0^*}$.

ASSUMPTION 3. In $\Theta \times \Gamma$, $M(\theta^*)$ is continuous and has a unique maximum at $\theta_0^*$.

REMARK 3. Under Assumptions 1 and 2, both Examples 1 and 2 satisfy Assumption 3. Actually in Example 1,

$$M(\theta^*) = -E_{\theta_0^*}[g(Y_i, \mathbf{x}(\theta^*, T_i))]$$

$$= -E_{\theta_0^*}[(Y_i - \mathbf{x}(\theta^*, T_i))^2] = -\sigma^2 - \int_0^T (\mathbf{x}(\theta^*, t) - \mathbf{x}(\theta_0^*, t))^2 Q(dt)$$

$$= M(\theta_0^*) - \int_0^T (\mathbf{x}(\theta^*, t) - \mathbf{x}(\theta_0^*, t))^2 Q(dt)$$

$$\leq M(\theta_0^*) - c \int_0^T (\mathbf{x}(\theta^*, t) - \mathbf{x}(\theta_0^*, t))^2 \, dt.$$

By Assumption 2, $\int_0^T (\mathbf{x}(\theta^*, t) - \mathbf{x}(\theta_0^*, t))^2 \, dt = 0$ if and only if $\theta^* = \theta_0^*$, hence $\theta_0^*$ is the unique maximizer of $M(\theta^*)$. For Example 2, the conditional probability

$$E_{\theta_0^*}[g(Y_i, \mathbf{x}(\theta^*, T_i))|T_i]$$

$$= -p(\theta_0^*, T_i) \log p(\theta^*, T_i) - (1 - p(\theta_0^*, T_i)) \log(1 - p(\theta^*, T_i)),$$

where

$$p(\theta^*, t) = \frac{e^{\mathbf{x}(\theta^*, t)}}{1 + e^{\mathbf{x}(\theta^*, t)}}.$$

Because for any fixed number $a \in (0, 1)$, the function

$$a \log x + (1 - a) \log (1 - x)$$

obtains its unique maximum at $a$ in $(0, 1)$,

$$M(\theta^*) = -E_{\theta_0^*}[g(Y_i, \mathbf{x}(\theta^*, T_i))]$$

$$= \int_0^T [p(\theta_0^*, t) \log p(\theta^*, t) + (1 - p(\theta_0^*, t)) \log (1 - p(\theta^*, t))] Q(dt)$$

$$= \int_0^T [p(\theta_0^*, t) \log p(\theta^*, t) + (1 - p(\theta_0^*, t)) \log (1 - p(\theta^*, t))] f(t) \, dt$$

obtains its maximum if and only if $\theta^* = \theta_0^*$.

ASSUMPTION 4. $g(y, x)$ is a nonnegative function and belongs to $C(\mathbb{R} \times \mathbb{R})$. If the $Y_i$ is not a bounded random variable, we assume that for any compact set $\Lambda \subset \mathbb{R}$,

$$(3.1) \qquad \liminf_{|y| \to \infty} \left[ \frac{1 + \inf_{x \in \Lambda} g(y, x)}{\sup_{x \in \Lambda} g(y, x)} \right] > 0.$$

REMARK 4. The latter statement means that for any two fixed points $x$, $x'$ in a compact subset, $g(y, x)$ and $g(y, x')$ are comparable as functions of $y$ when $|y| \to \infty$. Usually if $g(y, x)$ does not increase too fast as $|y| \to \infty$, the above condition is satisfied. For example, if $g(y, x)$ is a polynomial function of $y$ given $x$ and the coefficient of the highest order term is nonzero for all $x$, then (3.1) is true. Hence, Example 1 satisfies this assumption. In Example 2, $Y_i$ is bounded.

Given a compact subset $\Theta_0$ of $\Theta$ and a compact subset $\Gamma_0$ of $\Gamma$. Let

$$(3.2) \begin{aligned} r_n = \max \Bigg\{ & \sup_{\theta^* \in \Theta_0 \times \Gamma_0} \inf_{\mathbf{w} \in \mathbb{L}_n, \mathbf{w}(0) = x} \bigg[ \|\mathbf{x}(\theta^*, \cdot) - \mathbf{w}\|_\infty \\ & \qquad \vee \left\| \frac{d\mathbf{x}}{dt}(\theta^*, \cdot) - \frac{d\mathbf{w}}{dt} \right\|_\infty \bigg], \\ & \sup_{\theta^* \in \Theta_0 \times \Gamma_0} \inf_{\mathbf{v} \in \mathbb{L}_n, \mathbf{v}(0) = z} \bigg[ \|\mathbf{z}(\theta^*, \cdot) - \mathbf{v}\|_\infty \vee \left\| \frac{d\mathbf{z}}{dt}(\theta^*, \cdot) - \frac{d\mathbf{v}}{dt} \right\|_\infty \bigg] \Bigg\}, \end{aligned}$$

where $\theta^* = (\theta, x, z)$.

THEOREM 3.1. *Under Assumptions 1–4, suppose that $\lambda_n \to \infty$ and $r_n \to 0$ as $n \to \infty$. Then for any compact subset $\Theta_0$ of $\Theta$ and any compact subset $\Gamma_0$ of $\Gamma$,*

$$(3.3) \begin{aligned} & \sup_{\theta^* \in \Theta_0 \times \Gamma_0} \|\hat{\mathbf{x}}_n(\theta^*, \cdot) - \mathbf{x}(\theta^*, \cdot)\|_\infty \\ & \qquad \leq \left[ O_p \left( \frac{1}{\sqrt{\lambda_n}} \right) \sqrt{T} + 2T \sqrt{8(8K^2 + 2)} r_n \right] e^{2KT}, \end{aligned}$$

*where $K$ is a constant depending only on $\Theta_0 \times \Gamma_0$, $F$ and $G$.*

REMARK 5.

(1) There is an exponential function of $T$ in the above upper bound. This is because we use the penalty (1.4). According to the approximation theory for ODEs (see Antosiewicz [1]), for general ODEs, however small the approximations make such penalty (even we use the $L^\infty$-norms in the penalty), the difference between the approximations and the solutions could grow exponentially with $T$ increasing. Here is a simple example. Consider the equation,

$$\frac{dy}{dt} = y, \qquad y(0) = y_0.$$

The solution is $y(t) = y_0 e^t$. Now suppose we have an approximation $\hat{y}(t)$ to the solution satisfying

$$\frac{d\hat{y}}{dt}(t) - \hat{y}(t) = \varepsilon, \qquad t \geq 0,$$

where $\varepsilon$ is any fixed small number. Then we can get that

$$\hat{y}(t) = y_0 e^t + \varepsilon(e^t - 1), \qquad t \geq 0,$$

so

$$|\hat{y}(t) - y(t)| = |\varepsilon|(e^t - 1), \qquad t \geq 0.$$

The bound will increase very fast when $T$ is increasing. In the simulated data examples for FitzHugh–Nagumo equations in Ramsay et al. [21], they took $T = 20$. In this case, the above bound is too large to be useful for their sample size. However, the results in Ramsay et al. [21] and our simulation study indicate that when the smoothing parameter becomes large, the approximations to the solutions are very good. We will study this problem in the next section. Since $T$ is fixed, we can use this bound to get the asymptotic consistency and normality.

(2) If $\mathbb{L}_n$ is the space of cubic spline functions with knots $\tau^{(n)}$ by the proof of Lemma 1, we have

$$(3.4) \qquad\qquad r_n = O(|\tau^{(n)}|^3).$$

(3) The upper bound is the sum of two terms. The second term is a function of the distance between the approximation space and the true solutions. It does not depend on the sample and $\lambda$. This error term is due to the imperfect approximations of the basis expansion. The first term is due to using finite $\lambda$. If $\lambda$ is finite, the solutions to the penalized optimization problem (2.1) are affected by the sample, there are discrepancies between the solutions for finite $\lambda$ and the minimizer of $J$ (the solutions to the penalized optimization problem for $\lambda \to \infty$).

(4) In practice, numerical methods have to be used to calculate the penalties, which may lead to an increase in deviation of the approximations from the solutions. For example, Simpson's rule is used to compute the penalties in Ramsay et al. [21]. According to the error bound for Simpson's rule and the proof of Theorem 3.1, one more term should be added to the right-hand side of (3.3), $\tilde{K}|\tilde{\tau}|^{5/2}\sqrt{T}e^{2KT}$, where $\tilde{K}$ is a positive number depending on the derivatives of $\hat{\mathbf{x}}_n(\theta^*, t)$ up to order 5 and $\tilde{\tau}$ is the partition of $[0, T]$ for Simpson's rule. If the solutions are smooth enough and the knots for B-splines are included in the partition points for Simpson's method, the error can be well controlled by adding enough partition points.

In order to prove the consistency of the estimation $\hat{\theta}_n^*$, we replace Assumption 4 by a stronger assumption.

ASSUMPTION 5. $g(y, x)$ is a nonnegative function and belongs to $C^1(\mathbb{R} \times \mathbb{R})$ with

$$E_{\theta_0^*}\left[\left|\frac{\partial g}{\partial x}(Y_i, \mathbf{x}(\theta_0^*, T_i))\right|\right] < \infty.$$

If $Y_i$ is not a bounded random variables, we assume that for any compact set $\Lambda \subset \mathbb{R}$,

$$\liminf_{|y|\to\infty}\left[\frac{1 + \inf_{x\in\Lambda} g(y, x)}{\sup_{x\in\Lambda} g(y, x)}\right] > 0 \quad \text{and} \quad \liminf_{|y|\to\infty}\left[\frac{1 + \inf_{x\in\Lambda}|\partial g/\partial x(y, x)|}{\sup_{x\in\Lambda}|\partial g/\partial x(y, x)|}\right] > 0.$$

REMARK 6. Both Examples 1 and 2 satisfy this assumption.

THEOREM 3.2. *Suppose that Assumptions 1, 2, 3 and 5 hold and that $\hat{\theta}_n^*$ is uniformly tight. Suppose that $\lambda_n \to \infty$ and $r_n \to 0$ as $n \to \infty$, then $\hat{\theta}_n^*$ is consistent.*

REMARK 7.
(1) If $\mathbb{L}_n$ is the space of cubic spline functions with knots $\tau^{(n)}$. By (3.4), we only need $\lambda_n \to \infty$ and $|\tau^{(n)}| \to 0$ to obtain the consistency result.
(2) We say $\hat{\theta}_n^*$ is *uniformly tight* if for any $\varepsilon > 0$, there exists a compact set $\Theta_\varepsilon^* \subset \Theta \times \Gamma$ such that

$$\sup_n P(\hat{\theta}_n^* \notin \Theta_\varepsilon^*) < \varepsilon.$$

This is equivalent to stating that the probability of $\hat{\theta}_n$ going to the boundary or infinity is zero. The tightness assumption is essential for our proofs of the consistency and asymptotic normality. First, under this assumption, with a probability arbitrarily close to 1, the solutions can be uniformly bounded for any $n$. The weak convergence of the estimates also needs this assumption. If the parameter space $\Theta$ and the region $\Gamma$ where the initial values of ODEs are taken are bounded and closed, then $\hat{\theta}_n^*$ is automatically uniformly tight. For the general case, it is not easy to verify this assumption. Sometimes if $F, G$ go to infinity when $\theta$ goes to infinity or the boundary of $\Theta$, $\hat{\theta}_n^*$ is uniformly tight.

(3) The proofs of the consistency and the asymptotic normality are based on the results for $M$-estimators. Because $\hat{\theta}_n^*$ is the maximizer of $H_n(\hat{\mathbf{x}}_n(\theta^*, \cdot))$, it is natural to use $H_n(\hat{\mathbf{x}}_n(\theta^*, \cdot))$ as our criterion function. But $\hat{\mathbf{x}}_n$ is the solution of a penalized optimization problem, so it is an implicit function of $\theta^*$. The derivatives of $\hat{\mathbf{x}}_n$ with respect to $\theta^*$ have complicated forms and are difficult to analyze. Because we have obtained the upper bound for the difference between $\hat{\mathbf{x}}_n$ and $\mathbf{x}_n$, we can choose the criterion function $H_n(\mathbf{x}_n(\theta^*, \cdot))$ which can be more easily handled. Although $\hat{\theta}_n^*$ is not the maximizer of $H_n(\mathbf{x}_n(\theta^*, \cdot))$, we can control the difference between $H_n(\mathbf{x}_n(\hat{\theta}_n^*, \cdot))$ and $\max_{\theta^*} H_n(\mathbf{x}_n(\theta^*, \cdot))$. In other words, $\hat{\theta}_n^*$ nearly maximizes $H_n(\mathbf{x}_n(\theta^*, \cdot))$ (see Section 5.2 of van der Vaart [25]), so we can apply the results for $M$-estimators.

In order to prove the asymptotic normality of the estimator $\hat{\theta}_n^*$, we replace Assumptions 4 and 5 by a stronger assumption.

ASSUMPTION 6. $g(y, x)$ is a nonnegative function and belongs to $C^2(\mathbb{R} \times \mathbb{R})$ with

$$E_{\theta_0^*}\left[\left|\frac{\partial g}{\partial x}(Y_i, \mathbf{x}(\theta_0^*, T_i))\right|^2\right] < \infty \quad \text{and} \quad E_{\theta_0^*}\left[\left|\frac{\partial^2 g}{\partial x^2}(Y_i, \mathbf{x}(\theta_0^*, T_i))\right|\right] < \infty.$$

If the $Y_i$ is not a bounded random variable, we assume that for any compact set $\Lambda \subset \mathbb{R}$,

$$\liminf_{|y| \to \infty}\left[\frac{1 + \inf_{x \in \Lambda} g(y, x)}{\sup_{x \in \Lambda} g(y, x)}\right] > 0, \qquad \liminf_{|y| \to \infty}\left[\frac{1 + \inf_{x \in \Lambda} |\partial g/\partial x(y, x)|}{\sup_{x \in \Lambda} |\partial g/\partial x(y, x)|}\right] > 0$$

and

$$\liminf_{|y| \to \infty}\left[\frac{1 + \inf_{x \in \Lambda} |\partial^2 g/\partial x^2(y, x)|}{\sup_{x \in \Lambda} |\partial^2 g/\partial x^2(y, x)|}\right] > 0.$$

REMARK 8. Both Examples 1 and 2 satisfy this assumption.

Let the estimator $\mathring{\theta}_n^*$ be the maximizer of $H_n(\mathbf{x}(\theta^*, t))$. Note it is different from $\hat{\theta}_n^*$. If $g(y, \mathbf{x}(\theta, t))$ is the log density function of $(Y_i, T_i)$, then $\mathring{\theta}_n^*$ is the maximum likelihood estimator.

THEOREM 3.3. *Suppose that Assumptions 1, 2, 3 and 6 hold and that $\hat{\theta}_n^*$ and $\mathring{\theta}_n^*$ are uniformly tight. Suppose that*

$$\frac{\lambda_n}{n^2} \to \infty \quad \text{and} \quad r_n = o_p\left(\frac{1}{n}\right) \qquad \text{as } n \to \infty,$$

*and that the matrix*

$$V_{\theta_0^*} = -E_{\theta_0^*}\left[\frac{\partial g}{\partial x}(Y_i, \mathbf{x}(\theta_0^*, T_i))\, \frac{\partial^2 \mathbf{x}}{\partial \theta^* \partial \theta^{*T}}(\theta_0^*, T_i)\right.$$

$$\left. + \frac{\partial^2 g}{\partial x^2}(Y_i, \mathbf{x}(\theta_0^*, T_i))\, \frac{\partial \mathbf{x}}{\partial \theta^*}(\theta_0^*, T_i)\, \frac{\partial \mathbf{x}}{\partial \theta^*}(\theta_0^*, T_i)^T\right]$$

*is nonsingular. Then both $\sqrt{n}(\hat{\theta}_n^* - \theta_0^*)$ and $\sqrt{n}(\mathring{\theta}_n^* - \theta_0^*)$ are asymptotically normal with mean zero and the same asymptotic covariance matrix.*

REMARK 9.

(1) If $\mathbb{L}_n$ is the space of cubic spline functions with knots $\tau^{(n)}$. By (3.4), let

$$\frac{\lambda_n}{n^2} \to \infty \quad \text{and} \quad |\tau^{(n)}| = o_p\left(\frac{1}{n^{1/3}}\right) \qquad \text{as } n \to \infty.$$

Then the conditions on $\lambda_n$ and $r_n$ in Theorem 3.3 are satisfied.
(2) If $g(y, \mathbf{x}(\theta, t))$ is the log density function of $(Y_i, T_i)$, then $\mathring{\theta}_n^*$ is just the maximum likelihood estimation. Therefore, $\hat{\theta}_n^*$ is asymptotically efficient.
(3) The uniform tightness of $\mathring{\theta}_n^*$ is not needed in the proof of the asymptotic normality of $\hat{\theta}_n^*$.

## 4. The properties of the basis function approximations when $\lambda \to \infty$.
In this section, we study the finite sample behavior of the approximations of the solutions of (1.1) for a given sample and a given value of $\theta^*$. In Theorem 3.1, we provide a bound on the uniform norm of the difference between the approximations and the solutions. But this bound will grow exponentially with $T$ increasing due to the form of the penalty, which makes the bound useless in the finite-sample situation. It seems that the bound cannot be improved for general ODEs when the smoothing parameter is finite. However, the results in Ramsay et al. [21] and our simulation study indicate that when the smoothing parameter becomes large, the approximations to the solutions are quite good. Therefore, we let the smoothing parameter $\lambda$ go to infinity and study the limiting behavior of the approximations.

First we consider a simulated example. Consider the FitzHugh–Nagumo equations

$$\begin{aligned}
\frac{d\mathbf{V}}{dt}(t) &= c\left(\mathbf{V} - \frac{\mathbf{V}^3}{3} + \mathbf{R}\right), \\
\frac{d\mathbf{R}}{dt}(t) &= -\frac{1}{c}(\mathbf{V} - a + b\mathbf{R}).
\end{aligned} \tag{4.1}$$

The parameters in the system are $\theta = (a, b, c)$ and the time interval is $[0, 20]$. Suppose that the true parameter values are $\theta_0 = (0.2, 0.2, 3)$ and the estimates are $\hat{\theta} = (0.8, -0.5, 3.5)$. Let $(\mathbf{V}(\theta, \cdot), \mathbf{R}(\theta, \cdot))$ be the solutions of (4.1) with parameter $\theta$ and initial value $(1, -1)$. The data were simulated from the model

$$
\begin{aligned}
Y_{1i} &= \mathbf{V}(\theta_0, T_i) + \varepsilon_{1i}, \\
Y_{2i} &= \mathbf{R}(\theta_0, T_i) + \varepsilon_{2i},
\end{aligned}
\tag{4.2}
$$

where the samples were taken at times $0.0, 0.05, 0.10, \ldots, 20.0$ and $\varepsilon_{1i}, \varepsilon_{2i}$ were independent random variables with the same distributions $N(0, 0.5)$. We firstly plot the solutions for $\theta_0$ and $\hat{\theta}$, and the simulated data in Figure 1. Now, we fix the sample and let $\mathbb{L}$ be the cubic spline functions with knots at each sample point. Let the smoothing parameter $\lambda$ go to infinity. We plot the differences between the spline approximations $(\mathbf{X}_1(t), \mathbf{X}_2(t))$ and the solutions $(\mathbf{V}(t), \mathbf{R}(t))$ for $\hat{\theta} = (0.8, -0.5, 3.5)$ for different values of $\lambda$ in Figure 2. From this figure, we can see that the approximations to the solutions are quite good when $\lambda$ is large and there is no obvious time trend for the difference.

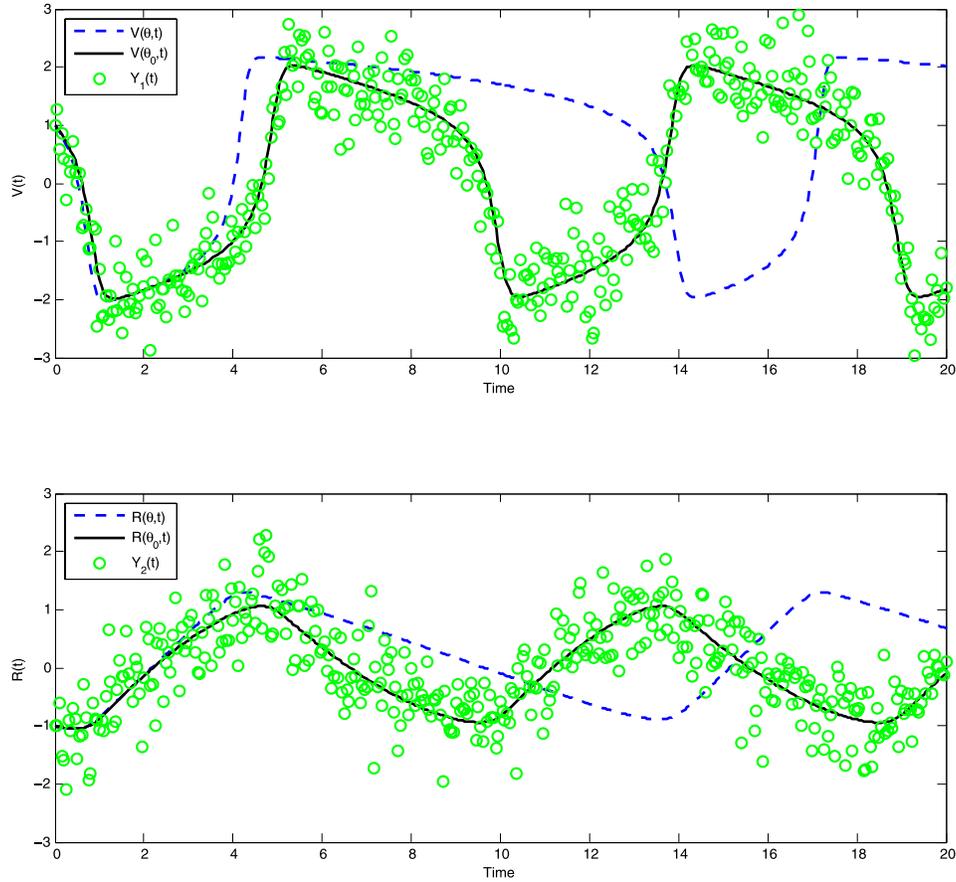

Fig. 1. *The solutions* $(\mathbf{V}, \mathbf{R})$ *of (4.1) for* $\theta_0 = (0.2, 0.2, 3)$ *and* $\theta = (0.8, -0.5, 3.5)$, *and the simulated data* $(Y_1, Y_2)$ *from (4.2).*

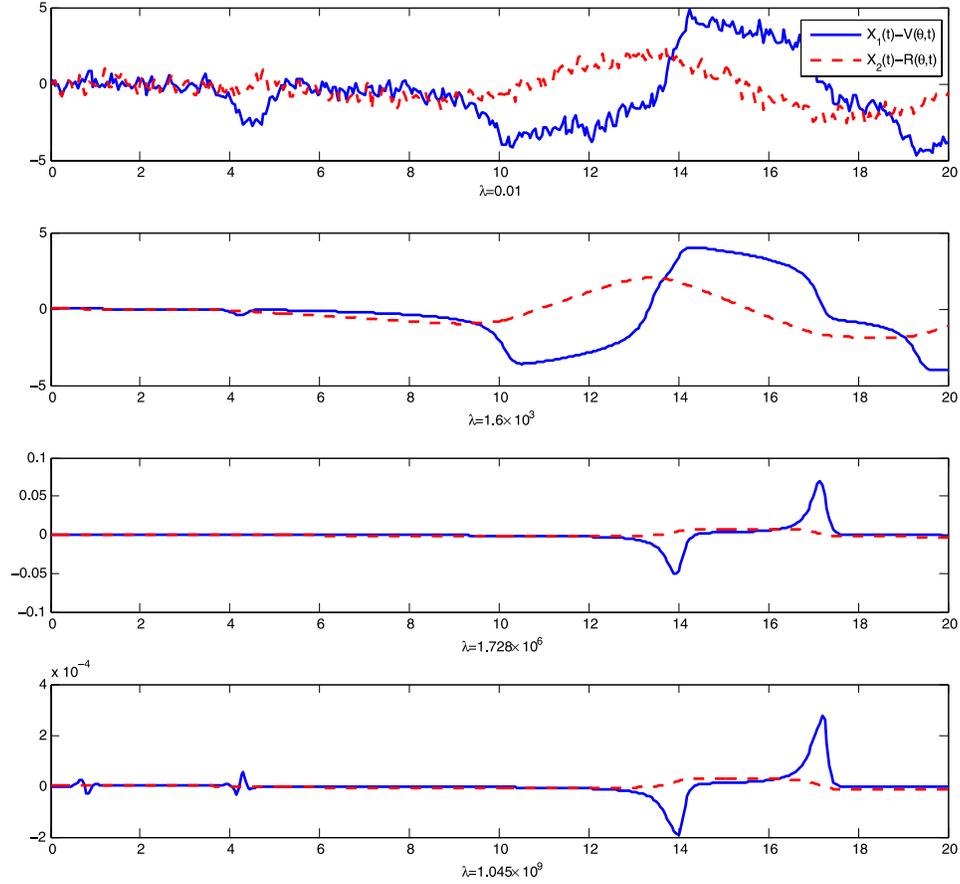

Fig. 2. *The differences between the spline approximations* $(\mathbf{X}_1(t), \mathbf{X}_2(t))$ *and the solution* $(\mathbf{V}(t), \mathbf{R}(t))$ *for* $\hat{\theta} = (0.8, -0.5, 3.5)$ *for different values of* $\lambda$. *The graphs have different y-axis scales.*

Fix the sample $(Y_1, T_1), \ldots, (Y_n, T_n)$ and the parameter $\theta^*$. Let $\mathbb{L}$ be a finite-dimensional linear subspace of $C^1[0, T]$. Define

$$r = \max\left\{ \inf_{\mathbf{w} \in \mathbb{L}, \mathbf{w}(0) = x} \left[ \|\mathbf{x}(\theta^*, \cdot) - \mathbf{w}\|_\infty \vee \left\| \frac{d\mathbf{x}}{dt}(\theta^*, \cdot) - \frac{d\mathbf{w}}{dt} \right\|_\infty \right],$$

(4.3)

$$\inf_{\mathbf{v}\in\mathbb{L},\mathbf{v}(0)=z}\left[\|\mathbf{z}(\theta^*,\cdot)-\mathbf{v}\|_\infty \vee \left\|\frac{d\mathbf{z}}{dt}(\theta^*,\cdot)-\frac{d\mathbf{v}}{dt}\right\|_\infty\right]\right\}.$$

Let $\lambda^{(m)}$ be a sequence of positive smoothing parameters which is strictly increasing and go to infinity. For each $m$, let

$$(\hat{\mathbf{x}}^{(m)},\hat{\mathbf{z}}^{(m)})\in \operatorname*{arg\,max}_{\substack{\hat{\mathbf{x}},\hat{\mathbf{z}}\in\mathbb{L}\\ \hat{\mathbf{x}}(0)=x,\hat{\mathbf{z}}(0)=z}} [H_n(\hat{\mathbf{x}})-\lambda^{(m)}J(\hat{\mathbf{x}},\hat{\mathbf{z}},\theta)].$$

Note that we suppress the subscript $n$ on $(\hat{\mathbf{x}}^{(m)}(\theta^*,t),\hat{\mathbf{z}}^{(m)}(\theta^*,t))$, $r$ and $\mathbb{L}$, which may depend on the sample, because the sample is fixed.

LEMMA 2. *For each $m$, we have*

$$H_n(\hat{\mathbf{x}}^{(m)})-H_n(\hat{\mathbf{x}}^{(m+1)})\le \lambda^{(m+1)}[J(\hat{\mathbf{x}}^{(m)},\hat{\mathbf{z}}^{(m)},\theta)-J(\hat{\mathbf{x}}^{(m+1)},\hat{\mathbf{z}}^{(m+1)},\theta)],$$

$$H_n(\hat{\mathbf{x}}^{(m)})-H_n(\hat{\mathbf{x}}^{(m+1)})\ge \lambda^{(m)}[J(\hat{\mathbf{x}}^{(m)},\hat{\mathbf{z}}^{(m)},\theta)-J(\hat{\mathbf{x}}^{(m+1)},\hat{\mathbf{z}}^{(m+1)},\theta)].$$

*Therefore, both $\{H_n(\hat{\mathbf{x}}^{(m)}):m\ge 1\}$ and $\{J(\hat{\mathbf{x}}^{(m)},\hat{\mathbf{z}}^{(m)},\theta):m\ge 1\}$ are decreasing sequences.*

PROOF. By the definitions of $(\hat{\mathbf{x}}^{(m)},\hat{\mathbf{z}}^{(m)})$ and $(\hat{\mathbf{x}}^{(m+1)},\hat{\mathbf{z}}^{(m+1)})$, we have

$$H_n(\hat{\mathbf{x}}^{(m)})-\lambda^{(m+1)}J(\hat{\mathbf{x}}^{(m)},\hat{\mathbf{z}}^{(m)},\theta)\le H_n(\hat{\mathbf{x}}^{(m+1)})-\lambda^{(m+1)}J(\hat{\mathbf{x}}^{(m+1)},\hat{\mathbf{z}}^{(m+1)},\theta),$$

$$H_n(\hat{\mathbf{x}}^{(m)})-\lambda^{(m)}J(\hat{\mathbf{x}}^{(m)},\hat{\mathbf{z}}^{(m)},\theta)\ge H_n(\hat{\mathbf{x}}^{(m+1)})-\lambda^{(m)}J(\hat{\mathbf{x}}^{(m+1)},\hat{\mathbf{z}}^{(m+1)},\theta).$$

The two inequalities in lemma follow immediately. Note that $\lambda^{(m+1)}>\lambda^{(m)}$, so we have

$$J(\hat{\mathbf{x}}^{(m)},\hat{\mathbf{z}}^{(m)},\theta)-J(\hat{\mathbf{x}}^{(m+1)},\hat{\mathbf{z}}^{(m+1)},\theta)\ge 0. \qquad \square$$

LEMMA 3.

$$\lim_{m\to\infty}J(\hat{\mathbf{x}}^{(m)},\hat{\mathbf{z}}^{(m)},\theta)=\inf_{\substack{\hat{\mathbf{x}},\hat{\mathbf{z}}\in\mathbb{L}\\ \hat{\mathbf{x}}(0)=x,\hat{\mathbf{z}}(0)=z}}J(\hat{\mathbf{x}},\hat{\mathbf{z}},\theta).$$

PROOF. By Lemma 2, $\{J(\hat{\mathbf{x}}^{(m)},\hat{\mathbf{z}}^{(m)},\theta):m\ge 1\}$ is a decreasing sequence and nonnegative, hence the limit exists. Now suppose that the equality is not true. Then we can find $(\tilde{\mathbf{x}},\tilde{\mathbf{z}})\in\mathbb{L}$ with $\tilde{\mathbf{x}}(0)=x,\tilde{\mathbf{z}}(0)=z$, and a number $\eta>0$ such that $\tilde{\mathbf{x}}(0)=x$

$$\lim_{m\to\infty}J(\hat{\mathbf{x}}^{(m)},\hat{\mathbf{z}}^{(m)},\theta)>J(\tilde{\mathbf{x}},\tilde{\mathbf{z}},\theta)+\eta.$$

For any $m$, we have

$$H_n(\hat{\mathbf{x}}^{(m)})-\lambda^{(m)}J(\hat{\mathbf{x}}^{(m)},\hat{\mathbf{z}}^{(m)},\theta)\ge H_n(\tilde{\mathbf{x}})-\lambda^{(m)}J(\tilde{\mathbf{x}},\tilde{\mathbf{z}},\theta).$$

Then

$$H_n(\hat{\mathbf{x}}^{(m)}) - H_n(\tilde{\mathbf{x}}) \geq \lambda^{(m)}[J(\hat{\mathbf{x}}^{(m)}, \hat{\mathbf{z}}^{(m)}, \theta) - J(\tilde{\mathbf{x}}, \tilde{\mathbf{z}}, \theta)] \geq \lambda^{(m)}\eta.$$

Note that $H_n(\cdot)$ is nonpositive, we have

$$-H_n(\tilde{\mathbf{x}}) \geq \lambda^{(m)}\eta.$$

The left-hand side is a fixed number and the right-hand side goes to infinity as $m \to \infty$. This is a contradiction. Hence, the lemma is true. $\square$

LEMMA 4.   *If both $\hat{\mathbf{x}}^{(m)}$ and $\hat{\mathbf{z}}^{(m)}$ are bounded sequences in $(\mathbb{L}, \|\cdot\|_\infty)$. Then for any subsequence $\{(\hat{\mathbf{x}}^{(m')}, \hat{\mathbf{z}}^{(m')})\} \subset \{(\hat{\mathbf{x}}^{(m)}, \hat{\mathbf{z}}^{(m)})\}$, there exist a further subsequence $\{(\hat{\mathbf{x}}^{(m'')}, \hat{\mathbf{z}}^{(m'')})\} \subset \{(\hat{\mathbf{x}}^{(m')}, \hat{\mathbf{z}}^{(m')})\}$ and*

$$(\hat{\mathbf{x}}^{(\infty)}, \hat{\mathbf{z}}^{(\infty)}) \in \underset{\substack{\hat{\mathbf{x}}, \hat{\mathbf{z}} \in \mathbb{L} \\ \hat{\mathbf{x}}(0) = x, \hat{\mathbf{z}} = z}}{\arg\min} \; J(\hat{\mathbf{x}}, \hat{\mathbf{z}}, \theta),$$

*such that*

$$\lim_{m'' \to \infty} \|\hat{\mathbf{x}}^{(m'')} - \hat{\mathbf{x}}^{(\infty)}\|_\infty = 0, \qquad \lim_{m'' \to \infty} \|\hat{\mathbf{z}}^{(m'')} - \hat{\mathbf{z}}^{(\infty)}\|_\infty = 0.$$

PROOF.   Because $(\mathbb{L}, \|\cdot\|_\infty)$ is a finite-dimensional subspace of $(C^1[0, T], \|\cdot\|_\infty)$, any bounded and closed subset of $(\mathbb{L}, \|\cdot\|_\infty)$ is compact. Then our conclusion follows from the fact that $J(\hat{\mathbf{x}}, \hat{\mathbf{z}}, \theta)$ is a continuous function of $(\hat{\mathbf{x}}, \hat{\mathbf{z}})$ in $(\mathbb{L}, \|\cdot\|_\infty)$ (note that it is not continuous in $(C^1[0, T], \|\cdot\|_\infty)$). $\square$

LEMMA 5.   *Suppose that the equations in (1.1) have unique solutions $(\mathbf{x}, \mathbf{z})$. For any $M > 0$ and $\delta > 0$, there exists a positive number $\varepsilon$ depending on $M$ and $\delta$, such that if $r < \varepsilon$, then for any*

$$(\hat{\mathbf{x}}^{(\infty)}, \hat{\mathbf{z}}^{(\infty)}) \in \underset{\substack{\hat{\mathbf{x}}, \hat{\mathbf{z}} \in \mathbb{L}, \hat{\mathbf{x}}(0) = x, \hat{\mathbf{z}} = z \\ \|\hat{\mathbf{x}}\|_\infty \leq M, \|\hat{\mathbf{z}}\|_\infty \leq M}}{\arg\min} \; J(\hat{\mathbf{x}}, \hat{\mathbf{z}}, \theta),$$

*we have*

$$\|\hat{\mathbf{x}}^{(\infty)} - \mathbf{x}\|_\infty < \delta, \qquad \|\hat{\mathbf{z}}^{(\infty)} - \mathbf{z}\|_\infty < \delta.$$

Now, we study the minimum points of $J$ in a neighborhood of the solutions. We can only give the result for the one-dimensional case and we need some properties of the basis functions. So we will assume in the next theorem that the subspace $\mathbb{L}$ is the space of B-spline functions with order at least 4. Let $\tau = (0 = t_1 < \cdots < t_{k_n} = T)$ be the knots of $\mathbb{L}$. Recall that

$$|\tau| = \max_{2 \leq i \leq k} |t_i - t_{i-1}|.$$

We define the mesh ratio

$$\kappa = \frac{\max_{2 \le i \le k} |t_i - t_{i-1}|}{\min_{2 \le i \le k} |t_i - t_{i-1}|} = \frac{|\tau|}{\min_{2 \le i \le k} |t_i - t_{i-1}|}.$$

Suppose that $\mathbf{x}$ is the solution of the equation

$$(4.4) \qquad \frac{d\mathbf{x}}{dt}(t) = F(\mathbf{x}(t), t), \qquad \mathbf{x}(0) = x,$$

where $F$ is a function on $\mathbb{R} \times \mathbb{R}$. We use $F_x, F_t$ to denote the partial derivatives of $F$ with respect to $x, t$, respectively. Define

$$J(\hat{\mathbf{x}}) = \int_0^T \left| \frac{d\hat{\mathbf{x}}}{dt}(t) - F(\hat{\mathbf{x}}(t), t) \right|^2 dt \qquad \forall \hat{\mathbf{x}} \in \mathbb{L}$$

and

$$(4.5) \qquad r = \inf_{\mathbf{w} \in \mathbb{L}, \mathbf{w}(0) = x} \left[ \|\mathbf{x} - \mathbf{w}\|_\infty \vee \left\| \frac{d\mathbf{x}}{dt} - \frac{d\mathbf{w}}{dt} \right\|_\infty \right].$$

THEOREM 4.1. *Assume that $\mathbf{x}$ is the unique solution of (4.4), $F$ has third-order continuous partial derivatives and*

$$F_x(\mathbf{x}(t), t) < 0 \qquad \forall 0 \le t \le T.$$

*Suppose that $\mathbb{L}$ is the cubic spline space with knots $\tau$ and $|\tau| \le 1$. Then there exists a positive number $\delta_2$ depending only on $\mathbf{x}$ and $F$. For any*

$$\hat{\mathbf{x}}_0 \in \operatorname*{arg\,min}_{\substack{\hat{\mathbf{x}} \in \mathbb{L}, \hat{\mathbf{x}}(0) = x \\ \|\hat{\mathbf{x}} - \mathbf{x}\|_\infty < \delta_2}} J(\hat{\mathbf{x}}),$$

*we have*

$$
\begin{aligned}
(4.6) \qquad \|\mathbf{x} - \hat{\mathbf{x}}_0\|_\infty < {}& \beta_1 \kappa \left\| \frac{d^2 \hat{\mathbf{x}}_0}{dt^2} \right\|_{L^2[0,T]} |\tau|^{1/2} + \beta_2 \kappa |\tau| \\
& + \kappa (4\sqrt{6\kappa} + \beta_3) \beta_4 \sqrt{T} |\tau|^{3/2} \\
& + \beta_5 \sqrt{T} |\tau|^3 + \beta_6 \sqrt{T} |\tau|^{7/2},
\end{aligned}
$$

*where $\beta_1, \beta_2, \beta_3, \beta_4, \beta_5, \beta_6$ are constants depending only on $\mathbf{x}$ and $F$.*

REMARK 10.

(1) Theorem 4.1 is true for any space of the B-spline functions with order larger than 4. However, the order of the bound on the right-hand side of (4.6) may be different for higher order B-spline functions. We conjecture that there are similar results for high-dimensional cases and more general equations.

(2) If $|\tau|$ is small enough, the bound is dominated by the first term which depends on $T$ only through $\|\frac{d^2\hat{x}_0}{dt^2}\|_{L^2[0,T]}$. If $\frac{d^2\hat{x}_0}{dt^2}$ is bounded by some fixed constant, then we can prove that the first term will be $O(|\tau|)$ and does not depend on $T$.

(3) Now, we can explain the pattern in Figure 2. According to Lemma 5, when $|\tau|$ is small enough, all the minimum points of $J$ will be in the neighborhood of the solution $(\mathbf{x}, \mathbf{z})$. By Theorem 4.1, all these minimum points will satisfy the bound (4.6). By Lemma 4, for any subsequence of $(\hat{\mathbf{x}}^{(m)}, \hat{\mathbf{z}}^{(m)})$, there exist a further subsequence converging to one of these minimum points. Hence, it is easy to show that $(\hat{\mathbf{x}}^{(m)}, \hat{\mathbf{z}}^{(m)})$ will satisfy the bound (4.6) for all large $m$.

(4) Here we outline the proof of this theorem. The first step is to derive a bound for $\|\mathbf{x}(t_0) - \hat{\mathbf{x}}_0(t_0)\|$ at any stationary point $t_0$ of the function $\mathbf{x}(t) - \hat{\mathbf{x}}_0(t)$, that is, any point $t_0$ such that $\frac{d\mathbf{x}}{dt}(t_0) - \frac{d\hat{\mathbf{x}}_0}{dt}(t_0) = 0$. In this step, we use the property of $\hat{\mathbf{x}}_0$ as a local minimizer of $J(\hat{\mathbf{x}})$. One of the key points in this step is based on the following observation: for any continuously differentiable function $f$ in $[0, T]$, if there is a zero point of $f$ between $a$ and $b$, then we have

$$\int_a^b f(t)\,dt = O(|b - a|^2).$$

We apply this observation to $\frac{d\mathbf{x}}{dt} - \frac{d\hat{\mathbf{x}}_0}{dt}$. In high-dimensional cases, the stationary points for different components of $\mathbf{x} - \hat{\mathbf{x}}_0$ may not be the same. We cannot extend the idea to high-dimensional cases. In this step, we just need the condition that $F_x(\mathbf{x}(t), t)$ is nonzero and do not require that it be strictly negative.

The next step is to give a bound for all the other points in $[0, T]$. In this step, we need the condition that $F_x(\mathbf{x}(t), t)$ is strictly negative to control the growth of $\|\mathbf{x}(t) - \hat{\mathbf{x}}_0(t)\|$ when $t$ gets away from a stationary point.

(5) In the proof, we use some properties of the B-spline bases. For example, the B-spline bases are stable bases and locally supported. Because we need to change the order of the differentiation and the integral in the proof, our proof cannot be applied to the penalty where $L^1$- or $L^\infty$-norm is used.

In Section 2, we prove the consistency and asymptotic normality based on the following idea, the likelihood functions can be well approximated if we have good approximations to the solutions of ODEs. Therefore, in practical applications, we should make the approximations close enough to the solutions. In some cases, if we do not want to solve the ODEs, we can use Theorem 4.1 to estimate the deviation of the spline approximations from solutions. We select a basis that is sufficiently rich to make the uniform norm in Theorem 4.1 very small. For any given $\lambda$ and $\theta^*$, if we want to estimate $\|\hat{\mathbf{x}}(\theta^*, \lambda, \cdot) - \mathbf{x}(\theta^*, \cdot)\|$, pick a sequence $\lambda_n \to \infty$. According to Lemma 4, we can find a subsequence $\lambda_{n'}$ and $\hat{\mathbf{x}}(\theta^*, \infty, \cdot)$ such that $\lim_{n' \to \infty} \hat{\mathbf{x}}(\theta^*, \lambda_{n'}, \cdot) = \hat{\mathbf{x}}(\theta^*, \infty, \cdot)$. By Theorem 4.1, we can use $\|\hat{\mathbf{x}}(\theta^*, \lambda, \cdot) - \hat{\mathbf{x}}(\theta^*, \infty, \cdot)\|$ to approximate $\|\hat{\mathbf{x}}(\theta^*, \lambda, \cdot) - \mathbf{x}(\theta^*, \cdot)\|$. In our simulation study, usually the sequence $\hat{\mathbf{x}}(\theta^*, \lambda_n, \cdot)$ converges, and there is no need to find the subsequence.

## 5. Proofs.

PROOF OF THEOREM 3.1. Given two compact sets $\Theta_0$ and $\Gamma_0$. Without loss of generality, we assume that $\Theta_0$ and $\Gamma_0$ are convex and contain $\theta_0$ and $(x_0, z_0)$, otherwise we can prove the conclusion for their convex hulls generated by $\Theta_0 \cup \{\theta_0\}$ and $\Gamma_0 \cup \{(x_0, z_0)\}$ which are still compact sets. Let $r_n$ be the number defined in (3.2). Since $r_n \to 0$, without loss of generality, we assume that $r_n \leq 1$. By the definition of $r_n$, for each $\theta^* \in \Theta_0 \times \Gamma_0$, there exist $\mathbf{w}_n(\theta^*, \cdot), \mathbf{v}_n(\theta^*, \cdot) \in \mathbb{L}_n$ such that

$$
\begin{aligned}
(5.1) \quad & \|\mathbf{x}(\theta^*, \cdot) - \mathbf{w}_n(\theta^*, \cdot)\|_\infty \vee \left\|\frac{d\mathbf{x}}{dt}(\theta^*, \cdot) - \frac{d\mathbf{w}_n}{dt}(\theta^*, \cdot)\right\|_\infty \leq 2r_n, \\
& \|\mathbf{z}(\theta^*, \cdot) - \mathbf{v}_n(\theta^*, \cdot)\|_\infty \vee \left\|\frac{d\mathbf{z}}{dt}(\theta^*, \cdot) - \frac{d\mathbf{v}_n}{dt}(\theta^*, \cdot)\right\|_\infty \leq 2r_n
\end{aligned}
$$

and $\mathbf{w}_n(\theta^*, 0) = \mathbf{x}(\theta^*, 0) = x$, $\mathbf{v}_n(\theta^*, 0) = \mathbf{z}(\theta^*, 0) = z$, where $\theta^* = (\theta, x, z)$. Because $(\mathbf{x}(\theta^*, t), \mathbf{z}(\theta^*, t))$ are continuous functions of $(\theta^*, t)$ and $\Theta_0 \times \Gamma_0 \times [0, T]$ is a compact set, there exists a positive number $R$ depending on $F, G$ and $\Theta_0, \Gamma_0$, such that

$$
(5.2) \quad |\mathbf{x}(\theta^*, t)| \leq R, \qquad |\mathbf{z}(\theta^*, t)| \leq R \qquad \forall (\theta^*, t) \in \Theta_0 \times \Gamma_0 \times [0, T].
$$

By (5.1),

$$
\begin{aligned}
(5.3) \quad |\mathbf{w}_n(\theta^*, t)| \leq R + 2r_n \leq R + 2, \qquad |\mathbf{v}_n(\theta^*, t)| \leq R + 2r_n \leq R + 2, \\
\forall (\theta^*, t) \in \Theta_0 \times \Gamma_0 \times [0, T].
\end{aligned}
$$

Since $F$, $G$ have continuous partial derivatives, we can find a positive number $K$ depending on $F, G$ and $\Theta_0, \Gamma_0$, such that

$$
\begin{aligned}
(5.4) \quad & |F(x, z, t, \theta) - F(x', z', t, \theta)| \leq K|x - x'| + K|z - z'|, \\
& |G(x, z, t, \theta) - G(x', z', t, \theta)| \leq K|x - x'| + K|z - z'| \\
& \forall |x| \leq R + 3, |x'| \leq R + 3, |z| \leq R + 3, |z'| \leq R + 3.
\end{aligned}
$$

We first prove a technical lemma.

LEMMA 6. $\sup_{\theta^* \in \Theta_0 \times \Gamma_0} \mathbb{P}_n g(Y, \mathbf{w}_n(\theta^*, \cdot)) = O_p(1)$.

PROOF. By (5.3), $\mathbf{w}_n(\theta^*, t)$ are uniformly bounded by $R + 2$ for all $(\theta^*, t) \in \Theta_0 \times \Gamma_0 \times [0, T]$. So if $Y_i$ is bounded, $\sup_{\theta^* \in \Theta_0 \times \Gamma_0} \mathbb{P}_n g(Y, \mathbf{w}_n(\theta^*, \cdot))$ is bounded since $g$ is continuous, and hence the lemma is true. Otherwise, by Assumption 4, we can find a positive number $\delta > 0$ such that

$$\delta \sup_{|x| \leq R+2} g(y, x) \leq \left(1 + \inf_{|x| \leq R+2} g(y, x)\right) \qquad \forall y \in \mathbb{R}.$$

Hence, we have

$$\sup_{\theta^* \in \Theta_0 \times \Gamma_0} \mathbb{P}_n g(Y, \mathbf{w}_n(\theta^*, \cdot)) \leq \frac{1}{\delta}(1 + \mathbb{P}_n g(Y, \mathbf{x}(\theta_0^*, \cdot)).$$

By the law of large numbers,

$$\mathbb{P}_n g(Y, \mathbf{x}(\theta_0^*, \cdot)) \to E_{\theta_0^*}[g(Y, \mathbf{x}(\theta_0^*, \cdot)] = M(\theta_0^*) < \infty,$$

$$\frac{1}{\delta}(1 + \mathbb{P}_n g(Y, \mathbf{x}(\theta_0^*, \cdot))) = O_p(1). \qquad \qquad \square$$

Now, we prove Theorem 3.1. For any $\theta^* \in \Theta_0 \times \Gamma_0$, by the definition of $(\hat{\mathbf{x}}_n(\theta^*, \cdot), \hat{\mathbf{z}}_n(\theta^*, \cdot))$, we have

$$H_n(\hat{\mathbf{x}}_n(\theta^*, \cdot)) - \lambda_n J(\hat{\mathbf{x}}_n(\theta^*, \cdot), \hat{\mathbf{z}}_n(\theta^*, \cdot), \theta)$$
$$\geq H_n(\mathbf{w}_n(\theta^*, \cdot)) - \lambda_n J(\mathbf{w}_n(\theta^*, \cdot), \mathbf{v}_n(\theta^*, \cdot), \theta).$$

Because $H_n(\hat{\mathbf{x}}_n(\theta^*, \cdot)) \leq 0$,

$$-\lambda_n J(\hat{\mathbf{x}}_n(\theta^*, \cdot), \hat{\mathbf{z}}_n(\theta^*, \cdot), \theta) \geq H_n(\mathbf{w}_n(\theta^*, \cdot)) - \lambda_n J(\mathbf{w}_n(\theta^*, \cdot), \mathbf{v}_n(\theta^*, \cdot), \theta).$$

Then

$$J(\hat{\mathbf{x}}_n(\theta^*, \cdot), \hat{\mathbf{z}}_n(\theta^*, \cdot), \theta)$$

$$\leq -\frac{1}{\lambda_n} H_n(\mathbf{w}_n(\theta^*, \cdot)) + J(\mathbf{w}_n(\theta^*, \cdot), \mathbf{v}_n(\theta^*, \cdot), \theta)$$

$$\leq \frac{1}{\lambda_n} \mathbb{P}_n g(Y, \mathbf{w}_n(\theta^*, \cdot)) + \left\| \frac{d\mathbf{w}_n}{dt}(\theta^*, \cdot) - F(\mathbf{w}_n(\theta^*, \cdot), \mathbf{v}_n(\theta^*, \cdot), t, \theta) \right\|_{L^2[0,T]}^2$$

$$\qquad + \left\| \frac{d\mathbf{v}_n}{dt}(\theta^*, \cdot) - G(\mathbf{w}_n(\theta^*, \cdot), \mathbf{v}_n(\theta^*, \cdot), t, \theta) \right\|_{L^2[0,T]}^2$$

$$= \frac{1}{\lambda_n} \mathbb{P}_n g(Y, \mathbf{w}_n(\theta^*, \cdot))$$

$$\qquad + \left\| \frac{d\mathbf{w}_n}{dt}(\theta^*, \cdot) - \frac{d\mathbf{x}}{dt}(\theta^*, \cdot) \right.$$

$$\qquad \left. + F(\mathbf{x}(\theta^*, \cdot), \mathbf{z}(\theta^*, \cdot), t, \theta) - F(\mathbf{w}_n(\theta^*, \cdot), \mathbf{v}_n(\theta^*, \cdot), t, \theta) \right\|_{L^2[0,T]}^2$$

$$+ \left\| \frac{d\mathbf{v}_n}{dt}(\theta^*,\cdot) - \frac{d\mathbf{z}}{dt}(\theta^*,\cdot) \right.$$

$$\left. + G(\mathbf{x}(\theta^*,\cdot),\mathbf{z}(\theta^*,\cdot),t,\theta) - G(\mathbf{w}_n(\theta^*,\cdot),\mathbf{v}_n(\theta^*,\cdot),t,\theta) \right\|^2_{L^2[0,T]}$$

$$(\mathbf{x},\mathbf{z} \text{ are solutions})$$

$$\leq \frac{1}{\lambda_n}\mathbb{P}_n g(Y,\mathbf{w}_n(\theta^*,\cdot)) + 2\left\| \frac{d\mathbf{w}_n}{dt}(\theta^*,\cdot) - \frac{d\mathbf{x}}{dt}(\theta^*,\cdot) \right\|^2_{L^2[0,T]}$$

$$+ 2\|F(\mathbf{x}(\theta^*,\cdot),\mathbf{z}(\theta^*,\cdot),t,\theta) - F(\mathbf{w}_n(\theta^*,\cdot),\mathbf{v}_n(\theta^*,\cdot),t,\theta)\|^2_{L^2[0,T]}$$

$$+ 2\left\| \frac{d\mathbf{v}_n}{dt}(\theta^*,\cdot) - \frac{d\mathbf{z}}{dt}(\theta^*,\cdot) \right\|^2_{L^2[0,T]}$$

$$+ 2\|G(\mathbf{x}(\theta^*,\cdot),\mathbf{z}(\theta^*,\cdot),t,\theta) - G(\mathbf{w}_n(\theta^*,\cdot),\mathbf{v}_n(\theta^*,\cdot),t,\theta)\|^2_{L^2[0,T]}$$

$$\leq \frac{1}{\lambda_n}\mathbb{P}_n g(Y,\mathbf{w}_n(\theta^*,\cdot)) + 2\left\| \frac{d\mathbf{w}_n}{dt}(\theta^*,\cdot) - \frac{d\mathbf{x}}{dt}(\theta^*,\cdot) \right\|^2_{L^2[0,T]}$$

$$+ 2\left\| \frac{d\mathbf{v}_n}{dt}(\theta^*,\cdot) - \frac{d\mathbf{z}}{dt}(\theta^*,\cdot) \right\|^2_{L^2[0,T]} + 8K^2\|\mathbf{w}_n(\theta^*,\cdot) - \mathbf{x}(\theta^*,\cdot)\|^2_{L^2[0,T]}$$

$$+ 8K^2\|\mathbf{v}_n(\theta^*,\cdot) - \mathbf{z}(\theta^*,\cdot)\|^2_{L^2[0,T]} \qquad \text{by (5.4)}$$

$$\leq \frac{1}{\lambda_n}\mathbb{P}_n g(Y,\mathbf{w}_n(\theta^*,\cdot)) + 2T\left\| \frac{d\mathbf{w}_n}{dt}(\theta^*,\cdot) - \frac{d\mathbf{x}}{dt}(\theta^*,\cdot) \right\|^2_\infty$$

$$+ 2T\left\| \frac{d\mathbf{v}_n}{dt}(\theta^*,\cdot) - \frac{d\mathbf{z}}{dt}(\theta^*,\cdot) \right\|^2_\infty + 8K^2 T\|\mathbf{w}_n(\theta^*,\cdot) - \mathbf{x}(\theta^*,\cdot)\|^2_\infty$$

$$+ 8K^2 T\|\mathbf{v}_n(\theta^*,\cdot) - \mathbf{z}(\theta^*,\cdot)\|^2_\infty \qquad \text{by (1.2)}$$

$$\leq \frac{1}{\lambda_n}\mathbb{P}_n g(Y,\mathbf{w}_n(\theta^*,\cdot))$$

$$+ (8K^2+2)T\left[ \left\| \frac{d\mathbf{w}_n}{dt}(\theta^*,\cdot) - \frac{d\mathbf{x}}{dt}(\theta^*,\cdot) \right\|_\infty \vee \|\mathbf{w}_n(\theta^*,\cdot) - \mathbf{x}(\theta^*,\cdot)\|_\infty \right]^2$$

$$+ (8K^2+2)T\left[ \left\| \frac{d\mathbf{v}_n}{dt}(\theta^*,\cdot) - \frac{d\mathbf{z}}{dt}(\theta^*,\cdot) \right\|_\infty \vee \|\mathbf{v}_n(\theta^*,\cdot) - \mathbf{z}(\theta^*,\cdot)\|_\infty \right]^2$$

$$\leq \frac{1}{\lambda_n}\mathbb{P}_n g(Y,\mathbf{w}_n(\theta^*,\cdot)) + 8T(8K^2+2)r_n^2 \qquad \text{by (5.1)}$$

and hence

$$\sup_{\theta^* \in \Theta_0 \times \Gamma_0} J(\hat{\mathbf{x}}_n(\theta^*, \cdot), \hat{\mathbf{z}}_n(\theta^*, \cdot), \theta)$$

$$= \frac{1}{\lambda_n} O_p(1) + 8T(8K^2 + 2)r_n^2$$

$$= O_p\left(\frac{1}{\lambda_n}\right) + 8T(8K^2 + 2)r_n^2 \qquad \text{by Lemma 6.}$$

By definition of $J$, we have

$$\sup_{\theta^* \in \Theta_0 \times \Gamma_0} \left\| \frac{d\hat{\mathbf{x}}_n}{dt}(\theta^*, \cdot) - F(\hat{\mathbf{x}}_n(\theta^*, \cdot), \hat{\mathbf{z}}_n(\theta^*, \cdot), t, \theta) \right\|_{L^2[0,T]}^2$$

$$\leq O_p\left(\frac{1}{\lambda_n}\right) + 8T(8K^2 + 2)r_n^2,$$

$$\sup_{\theta^* \in \Theta_0 \times \Gamma_0} \left\| \frac{d\hat{\mathbf{z}}_n}{dt}(\theta^*, \cdot) - G(\hat{\mathbf{x}}_n(\theta^*, \cdot), \hat{\mathbf{z}}_n(\theta^*, \cdot), t, \theta) \right\|_{L^2[0,T]}^2$$

$$\leq O_p\left(\frac{1}{\lambda_n}\right) + 8T(8K^2 + 2)r_n^2.$$

Therefore, by (1.3),

$$\sup_{\theta^* \in \Theta_0 \times \Gamma_0} \left\| \hat{\mathbf{x}}_n(\theta^*, t) - x - \int_0^t F(\hat{\mathbf{x}}_n(\theta^*, s), \hat{\mathbf{z}}_n(\theta^*, s), s, \theta)\, ds \right\|_\infty$$

(5.5)
$$\leq \sqrt{O_p\left(\frac{1}{\lambda_n}\right)T + 8T^2(8K^2 + 2)r_n^2}$$

$$\leq O_p\left(\frac{1}{\sqrt{\lambda_n}}\right)\sqrt{T} + T\sqrt{8(8K^2 + 2)}r_n,$$

and similarly,

$$\sup_{\theta^* \in \Theta_0 \times \Gamma_0} \left\| \hat{\mathbf{z}}_n(\theta^*, t) - z - \int_0^t G(\hat{\mathbf{x}}_n(\theta^*, s), \hat{\mathbf{z}}_n(\theta^*, s), s, \theta)\, ds \right\|_\infty$$

(5.6)
$$\leq O_p\left(\frac{1}{\sqrt{\lambda_n}}\right)\sqrt{T} + T\sqrt{8(8K^2 + 2)}r_n.$$

Define

$$A_n(\theta^*, t) = \hat{\mathbf{x}}_n(\theta^*, t) - x - \int_0^t F(\hat{\mathbf{x}}_n(\theta^*, s), \hat{\mathbf{z}}_n(\theta^*, s), s, \theta)\, ds,$$

(5.7)
$$B_n(\theta^*, t) = \hat{\mathbf{z}}_n(\theta^*, t) - z - \int_0^t G(\hat{\mathbf{x}}_n(\theta^*, s), \hat{\mathbf{z}}_n(\theta^*, s), s, \theta)\, ds.$$

We will show that if $\sup_{\theta^* \in \Theta_0 \times \Gamma_0} \|A_n(\theta^*, \cdot)\|_\infty$ and $\sup_{\theta^* \in \Theta_0 \times \Gamma_0} \|B_n(\theta^*, \cdot)\|_\infty$ are small enough, then $|\hat{\mathbf{x}}_n(\theta^*, t)| \leq R + 3$ and $|\hat{\mathbf{z}}_n(\theta^*, t)| \leq R + 3$ for all $t \in [0, T]$, $\theta^* \in \Theta_0 \times \Gamma_0$. Because $(\mathbf{x}(\theta^*, \cdot), \mathbf{z}(\theta^*, \cdot))$ are solutions of (1.1) with initial values $(x, z)$, by integrating the equations (1.1), we have

$$\mathbf{x}(\theta^*, t) - x - \int_0^t F(\mathbf{x}(\theta^*, s), \mathbf{z}(\theta^*, s), s, \theta) \, ds = 0,$$

$$\mathbf{z}(\theta^*, t) - z - \int_0^t G(\mathbf{x}(\theta^*, s), \mathbf{z}(\theta^*, s), s, \theta) \, ds = 0.$$

Then we subtract them from (5.7),

$$\begin{aligned} A_n(\theta^*, t) = {} & [\hat{\mathbf{x}}_n(\theta^*, t) - \mathbf{x}(\theta^*, t)] \\ & - \int_0^t [F(\hat{\mathbf{x}}_n(\theta^*, s), \hat{\mathbf{z}}_n(\theta^*, s), s, \theta) - F(\mathbf{x}(\theta^*, s), \mathbf{z}(\theta^*, s), s, \theta)] \, ds, \end{aligned}$$

$$\begin{aligned} B_n(\theta^*, t) = {} & [\hat{\mathbf{z}}_n(\theta^*, t) - \mathbf{z}(\theta^*, t)] \\ & - \int_0^t [G(\hat{\mathbf{x}}_n(\theta^*, s), \hat{\mathbf{z}}_n(\theta^*, s), s, \theta) - G(\mathbf{x}(\theta^*, s), \mathbf{z}(\theta^*, s), s, \theta)] \, ds. \end{aligned}$$

So

$$\begin{aligned} & |\hat{\mathbf{x}}_n(\theta^*, t) - \mathbf{x}(\theta^*, t)| \\ & \quad \leq \int_0^t |F(\hat{\mathbf{x}}_n(\theta^*, s), \hat{\mathbf{z}}_n(\theta^*, s), s, \theta) - F(\mathbf{x}(\theta^*, s), \mathbf{z}(\theta^*, s), s, \theta)| \, ds \\ & \qquad + |A_n(\theta^*, t)|, \\ & |\hat{\mathbf{z}}_n(\theta^*, t) - \mathbf{z}(\theta^*, t)| \\ & \quad \leq \int_0^t |G(\hat{\mathbf{x}}_n(\theta^*, s), \hat{\mathbf{z}}_n(\theta^*, s), s, \theta) - G(\mathbf{x}(\theta^*, s), \mathbf{z}(\theta^*, s), s, \theta)| \, ds \\ & \qquad + |B_n(\theta^*, t)|. \end{aligned}$$

Define $\tau_{\theta^*} = [\inf\{t \geq 0, |\hat{\mathbf{x}}_n(\theta^*, t)| \geq R + 3, \text{ or } |\hat{\mathbf{z}}_n(\theta^*, t)| \geq R + 3\}] \wedge T$, where for any $a, b \in \mathbb{R}$, $a \wedge b = \min(a, b)$. Then for any $0 \leq s \leq \tau_{\theta^*}$, and $\theta^* \in \Theta_0 \times \Gamma_0$, $|\hat{\mathbf{x}}_n(\theta^*, s)| \leq R + 3$, $|\hat{\mathbf{z}}_n(\theta^*, s)| \leq R + 3$. By applying (5.4), we have for any $0 \leq t \leq \tau_{\theta^*}$,

$$\begin{aligned} & |\hat{\mathbf{x}}_n(\theta^*, t) - \mathbf{x}(\theta^*, t)| \\ & \quad \leq \int_0^t |F(\hat{\mathbf{x}}_n(\theta^*, s), \hat{\mathbf{z}}_n(\theta^*, s), s, \theta) - F(\mathbf{x}(\theta^*, s), \mathbf{z}(\theta^*, s), s, \theta)| \, ds \\ & \qquad + |A_n(\theta^*, t)| \\ & \quad \leq \int_0^t |F(\hat{\mathbf{x}}_n(\theta^*, s), \hat{\mathbf{z}}_n(\theta^*, s), s, \theta) - F(\mathbf{x}(\theta^*, s), \mathbf{z}(\theta^*, s), s, \theta)| \, ds \end{aligned}$$

$$+ \sup_{\theta^* \in \Theta_0 \times \Gamma_0} \|A_n(\theta^*, \cdot)\|_\infty$$

$$\leq K \int_0^t |\hat{\mathbf{x}}_n(\theta^*, s) - \mathbf{x}(\theta^*, s)| \, ds$$

$$+ K \int_0^t |\hat{\mathbf{z}}_n(\theta^*, s) - \mathbf{z}(\theta^*, s)| \, ds$$

$$+ \sup_{\theta^* \in \Theta_0 \times \Gamma_0} \|A_n(\theta^*, \cdot)\|_\infty.$$

Similarly,

$$|\hat{\mathbf{z}}_n(\theta^*, t) - \mathbf{z}(\theta^*, t)|$$

$$\leq K \int_0^t |\hat{\mathbf{x}}_n(\theta^*, s) - \mathbf{x}(\theta^*, s)| \, ds + K \int_0^t |\hat{\mathbf{z}}_n(\theta^*, s) - \mathbf{z}(\theta^*, s)| \, ds$$

$$+ \sup_{\theta^* \in \Theta_0 \times \Gamma_0} \|B_n(\theta^*, \cdot)\|_\infty.$$

We add the above two inequalities together

$$|\hat{\mathbf{x}}_n(\theta^*, t) - \mathbf{x}(\theta^*, t)| + |\hat{\mathbf{z}}_n(\theta^*, t) - \mathbf{z}(\theta^*, t)|$$

$$\leq 2K \int_0^t [|\hat{\mathbf{x}}_n(\theta^*, s) - \mathbf{x}(\theta^*, s)| + |\hat{\mathbf{z}}_n(\theta^*, s) - \mathbf{z}(\theta^*, s)|] \, ds$$

$$+ \sup_{\theta^* \in \Theta_0 \times \Gamma_0} \|A_n(\theta^*, \cdot)\|_\infty + \sup_{\theta^* \in \Theta_0 \times \Gamma_0} \|B_n(\theta^*, \cdot)\|_\infty.$$

It follows from Gronwall's inequality that

$$|\hat{\mathbf{x}}_n(\theta^*, t) - \mathbf{x}(\theta^*, t)| + |\hat{\mathbf{z}}_n(\theta^*, t) - \mathbf{z}(\theta^*, t)|$$

(5.8)
$$\leq \left[ \sup_{\theta^* \in \Theta_0 \times \Gamma_0} \|A_n(\theta^*, \cdot)\|_\infty + \sup_{\theta^* \in \Theta_0 \times \Gamma_0} \|B_n(\theta^*, \cdot)\|_\infty \right] e^{2K\tau_{\theta^*}}$$

$$\leq \left[ \sup_{\theta^* \in \Theta_0 \times \Gamma_0} \|A_n(\theta^*, \cdot)\|_\infty + \sup_{\theta^* \in \Theta_0 \times \Gamma_0} \|B_n(\theta^*, \cdot)\|_\infty \right] e^{2KT}$$

$$\forall 0 \leq t \leq \tau_{\theta^*}, \theta^* \in \Theta_0 \times \Gamma_0.$$

By (5.5), (5.6) and (5.7), as $n \to \infty$, we have

$$\sup_{\theta^* \in \Theta_0 \times \Gamma_0} \|A_n(\theta^*, \cdot)\|_\infty + \sup_{\theta^* \in \Theta_0 \times \Gamma_0} \|B_n(\theta^*, \cdot)\|_\infty \to 0.$$

Hence, when $n$ is large enough, we have

$$|\hat{\mathbf{x}}_n(\theta^*, t) - \mathbf{x}(\theta^*, t)| + |\hat{\mathbf{z}}_n(\theta^*, t) - \mathbf{z}(\theta^*, t)| \leq 1.$$

By (5.2), $|\hat{\mathbf{x}}_n(\theta^*, t)| \leq |\mathbf{x}(\theta^*, t)| + 1 \leq R + 1 < R + 3$ and $|\hat{\mathbf{z}}_n(\theta^*, t)| \leq |\mathbf{z}(\theta^*, t)| + 1 \leq R + 1 < R + 3$ for all $t \in [0, \tau_{\theta^*}], \theta^* \in \Theta_0 \times \Gamma_0$. By the definition of $\tau_{\theta^*}$,

we must have $\tau_{\theta^*} = T$ for all $\theta^* \in \Theta_0 \times \Gamma_0$. By (5.5), (5.6), (5.7) and (5.8), we have

$$\sup_{\theta^* \in \Theta_0 \times \Gamma_0} \|\hat{\mathbf{x}}_n(\theta^*, \cdot) - \mathbf{x}(\theta^*, \cdot)\|_\infty \leq \left[ O_p\left(\frac{1}{\sqrt{\lambda_n}}\right)\sqrt{T} + 2T\sqrt{8(8K^2+2)}r_n \right] e^{2KT}.$$

$\square$

PROOF OF THEOREM 3.2. For any two given compact sets $\Theta_0 \subset \Theta$ and $\Gamma_0 \subset \Gamma$, first, we show

$$(5.9) \quad \begin{aligned} \sup_{\theta^* \in \Theta_0 \times \Gamma_0} & |H_n(\hat{\mathbf{x}}_n(\theta^*, \cdot)) - H_n(\mathbf{x}(\theta^*, \cdot))| \\ &= \left[ O_p\left(\frac{1}{\sqrt{\lambda_n}}\right)\sqrt{T} + 2T\sqrt{8(8K^2+2)}r_n \right] e^{2KT} O_p(1) + o_p\left(\frac{1}{n}\right). \end{aligned}$$

Second, define

$$M_n(\theta^*) = H_n(\mathbf{x}_n(\theta^*, \cdot)) = -\frac{1}{n}\sum_{i=1}^n g(Y_i, \mathbf{x}(\theta^*, T_i)), \qquad \theta^* \in \Theta \times \Gamma.$$

We show that

$$(5.10) \quad \sup_{\theta^* \in \Theta_0 \times \Gamma_0} |M_n(\theta^*) - M(\theta^*)| = o_p(1).$$

Note that $M(\theta^*) = -E_{\theta_0^*}[g(Y_i, \mathbf{x}(\theta^*, T_i))]$. Finally, we show that

$$\hat{\theta}_n^* \to \theta_0^* \qquad \text{in probability.}$$

Now let us prove (5.9) and (5.10). Without loss of generality, we assume that $\Theta_0$ and $\Gamma_0$ are convex and contain $\theta_0$ and $(x_0, z_0)$. According to Assumption 2, $(\mathbf{x}(\theta^*, t), \mathbf{z}(\theta^*, t))$ are continuous functions of $(\theta^*, t)$. Because $\Theta_0 \times \Gamma_0 \times [0, T]$ is a compact set, there exists a positive number $R$ depending on $\Theta_0$ and $\Gamma_0$, such that

$$|\mathbf{x}(\theta^*, t)| \leq R, \qquad |\mathbf{z}(\theta^*, t)| \leq R \qquad \forall (\theta^*, t) \in \Theta_0 \times \Gamma_0 \times [0, T].$$

Define $V_n = \sup_{\theta^* \in \Theta_0 \times \Gamma_0} \|\hat{\mathbf{x}}_n(\theta^*, \cdot) - \mathbf{x}(\theta^*, \cdot)\|_\infty$. From Theorem 3.1,

$$(5.11) \quad V_n \leq \left[ O_p\left(\frac{1}{\sqrt{\lambda_n}}\right)\sqrt{T} + 2T\sqrt{8(8K^2+2)}r_n \right] e^{2KT} = o_p(1).$$

By Assumption 5, there exists a positive number $\delta$, such that

$$(5.12) \quad \sup_{|x| \leq R+1} \left| \frac{\partial g}{\partial x}(y, x) \right| \leq \frac{1}{\delta}\left( 1 + \inf_{|x| \leq R+1} \left| \frac{\partial g}{\partial x}(y, x) \right| \right) \qquad \forall y \in \mathbb{R}.$$

Then we have

$$|H_n(\hat{\mathbf{x}}_n(\theta^*, \cdot)) - H_n(\mathbf{x}(\theta^*, \cdot))|$$

$$\leq \frac{1}{n} \sum_{i=1}^{n} |g(Y_i, \hat{\mathbf{x}}_n(\theta^*, T_i)) - g(Y_i, \mathbf{x}(\theta^*, T_i))|$$

$$\leq \frac{1}{n} \sum_{i=1}^{n} |g(Y_i, \hat{\mathbf{x}}_n(\theta^*, T_i)) - g(Y_i, \mathbf{x}(\theta^*, T_i))| \mathbf{1}_{[V_n \leq 1]}$$

$$+ \frac{1}{n} \sum_{i=1}^{n} |g(Y_i, \hat{\mathbf{x}}_n(\theta^*, T_i)) - g(Y_i, \mathbf{x}(\theta^*, T_i))| \mathbf{1}_{[V_n > 1]}$$

$$\leq \frac{1}{n} \sum_{i=1}^{n} \left[ \sup_{|x| \leq R+1} \left| \frac{\partial g}{\partial x}(Y_i, x) \right| \right] \|\hat{\mathbf{x}}_n(\theta^*, \cdot) - \mathbf{x}(\theta^*, \cdot)\|_\infty \mathbf{1}_{[V_n \leq 1]}$$

$$+ \frac{1}{n} \sum_{i=1}^{n} |g(Y_i, \hat{\mathbf{x}}_n(\theta^*, T_i)) - g(Y_i, \mathbf{x}(\theta^*, T_i))| \mathbf{1}_{[V_n > 1]}$$

$$\leq \frac{1}{n} \sum_{i=1}^{n} \frac{1}{\delta} \left[ 1 + \inf_{|x| \leq R+1} \left| \frac{\partial g}{\partial x}(Y_i, x) \right| \right] V_n \mathbf{1}_{[V_n \leq 1]} \qquad \text{by (5.12)}$$

$$+ \frac{1}{n} \sum_{i=1}^{n} |g(Y_i, \hat{\mathbf{x}}_n(\theta^*, T_i)) - g(Y_i, \mathbf{x}(\theta^*, T_i))| \mathbf{1}_{[V_n > 1]}$$

$$\leq \frac{1}{n} \sum_{i=1}^{n} \frac{1}{\delta} \left[ 1 + \left| \frac{\partial g}{\partial x}(Y_i, \mathbf{x}(\theta_0^*, T_i)) \right| \right] V_n$$

$$+ \frac{1}{n} \sum_{i=1}^{n} |g(Y_i, \hat{\mathbf{x}}_n(\theta^*, T_i)) - g(Y_i, \mathbf{x}(\theta^*, T_i))| \mathbf{1}_{[V_n > 1]}.$$

So

$$\sup_{\theta^* \in \Theta_0 \times \Gamma_0} |H_n(\hat{\mathbf{x}}_n(\theta^*, \cdot)) - H_n(\mathbf{x}(\theta^*, \cdot))|$$

$$(5.13) \qquad \leq \frac{1}{n} \sum_{i=1}^{n} \frac{1}{\delta} \left[ 1 + \left| \frac{\partial g}{\partial x}(Y_i, \mathbf{x}(\theta_0^*, T_i)) \right| \right] V_n$$

$$+ \sup_{\theta^* \in \Theta_0 \times \Gamma_0} \left[ \frac{1}{n} \sum_{i=1}^{n} |g(Y_i, \hat{\mathbf{x}}_n(\theta^*, T_i)) - g(Y_i, \mathbf{x}(\theta^*, T_i))| \right] \mathbf{1}_{[V_n > 1]}.$$

By the law of large numbers and Assumption 5,

$$\frac{1}{n}\sum_{i=1}^{n}\left[1+\left|\frac{\partial g}{\partial x}(Y_i,\mathbf{x}(\theta_0^*,T_i))\right|\right]\to 1+E_{\theta_0^*}\left|\frac{\partial g}{\partial x}(Y_i,\mathbf{x}(\theta_0^*,T_i))\right|<\infty$$

in probability, so it is $O_p(1)$. From (5.11), the second term on the right-hand side of (5.13) is not zero only in the event $[V_n>1]$ whose probability goes to zero, so it is $o_p(\frac{1}{n})$. Equation (5.9) has been proved. The equality (5.10) follows from the lemma below.

LEMMA 7. *The two classes $\{\mathbf{x}(\theta^*,\cdot),\theta^*\in\Theta_0\times\Gamma_0\}$ and $\{g(\cdot,\mathbf{x}(\theta^*,\cdot)),\theta^*\in\Theta_0\times\Gamma_0\}$ are both $P_{\theta_0^*}$-Glivenko–Cantelli.*

PROOF. For any $\theta^{*\prime},\theta^{*\prime\prime}\in\Theta_0\times\Gamma_0$, let $\theta^{*\prime}=(\theta',x',z')$ and $\theta^{*\prime\prime}=(\theta'',x'',z'')$. Since $\Theta_0\times\Gamma_0$ is convex, by Taylor expansion, we have

$$|\mathbf{x}(\theta^{*\prime},t)-\mathbf{x}(\theta^{*\prime\prime},t)|\leq\left[\sup_{\theta^*\in\Theta_0\times\Gamma_0,0\leq s\leq T}\left|\frac{\partial\mathbf{x}}{\partial\theta}(\theta^*,s)\right|\right]|\theta'-\theta''|$$

$$+\left[\sup_{\theta^*\in\Theta_0\times\Gamma_0,0\leq s\leq T}\left|\frac{\partial\mathbf{x}}{\partial x}(\theta^*,s)\right|\right]|x'-x''|$$

$$+\left[\sup_{\theta^*\in\Theta_0\times\Gamma_0,0\leq s\leq T}\left|\frac{\partial\mathbf{x}}{\partial z}(\theta^*,s)\right|\right]|z'-z''|.$$

By Assumption 2, $\frac{\partial\mathbf{x}}{\partial\theta}$, $\frac{\partial\mathbf{x}}{\partial x}$, $\frac{\partial\mathbf{x}}{\partial z}$ are all continuous functions of $(\theta^*,t)$. Therefore, they are all bounded in the compact set $\Theta_0\times\Gamma_0\times[0,T]$. We can find a positive constant $C$, such that

$$|\mathbf{x}(\theta^{*\prime},t)-\mathbf{x}(\theta^{*\prime\prime},t)|\leq C|\theta^{*\prime}-\theta^{*\prime\prime}|$$

$$\forall\theta^*\in\Theta_0\times\Gamma_0,0\leq t\leq T.$$

That is, the class $\{\mathbf{x}(\theta^*,\cdot),\theta^*\in\Theta_0\times\Gamma_0\}$ is Lipschitz in $\Theta_0\times\Gamma_0$. It follows from the Theorem 2.7.11 in the book of van der Vaart and Wellner [26] that the $L_1(P_{\theta_0^*})$-bracketing number is bounded by the covering number $N(\varepsilon,\Theta_0\times\Gamma_0,|\cdot|)$ of $\Theta_0\times\Gamma_0$, where $|\cdot|$ is the Euclidean distance. Because $\Theta_0\times\Gamma_0$ is a bounded subset in $\mathbb{R}^{d+2}$,

$$N(\varepsilon,\Theta_0\times\Gamma_0,|\cdot|)\leq\text{constant}\times\left(\frac{1}{\varepsilon}\right)^{d+2}.$$

Then our lemma follows from the Theorem 2.4.1 in the book of van der Vaart and Wellner [26]. Similar to the derivation of (5.13), we can get

$$|g(y,\mathbf{x}(\theta^{*\prime},t))-g(y,\mathbf{x}(\theta^{*\prime\prime},t))|$$

$$(5.14) \qquad \leq \frac{1}{\delta} \left[ 1 + \left| \frac{\partial g}{\partial x}(y, \mathbf{x}(\theta_0^*, t)) \right| \right] |\mathbf{x}(\theta^{*\prime}, t) - \mathbf{x}(\theta^{*\prime\prime}, t)|$$

$$\leq \frac{C}{\delta} \left[ 1 + \left| \frac{\partial g}{\partial x}(y, \mathbf{x}(\theta_0^*, t)) \right| \right] |\theta^{*\prime} - \theta^{*\prime\prime}|.$$

By Assumption 5, $|\frac{\partial g}{\partial x}(y, \mathbf{x}(\theta_0^*, t))|$ has a finite expectation. Hence, by the same argument for $\mathbf{x}(\theta^*, \cdot)$, we can get the conclusion for $g(\cdot, \mathbf{x}(\theta^*, \cdot))$.  $\square$

Because $\hat{\theta}_n^*$ is uniformly tight, for any $\varepsilon > 0$, there exist compact sets $\Theta_0$ and $\Gamma_0$ such that

$$(5.15) \qquad P_{\theta_0^*}(\Pi) \geq 1 - \varepsilon,$$

where

$$(5.16) \qquad \Pi = [\hat{\theta}_n^* \in \Theta_0 \times \Gamma_0 \text{ for all } n].$$

Without loss of generality, we assume that $\Theta_0$ and $\Gamma_0$ contain $\theta_0$ and $(x_0, z_0)$. Let $\eta$ be any positive number. By Assumption 3, $\theta_0^*$ is the unique maximum point of $M(\theta^*)$, hence

$$(5.17) \qquad \gamma = M(\theta_0^*) - \sup_{\theta^* \in \Theta_0 \times \Gamma_0, |\theta^* - \theta_0^*| \geq \eta} M(\theta^*) > 0$$

due to the compactness of $\Theta_0 \times \Gamma_0$ and the continuity of $M$. Because for $\hat{\theta}_n^* \in \Theta_0 \times \Gamma_0$,

$$M(\theta_0^*) - M(\hat{\theta}_n^*)$$
$$= M_n(\theta_0^*) - M_n(\hat{\theta}_n^*) + o_p(1) \qquad \text{by } (5.10)$$
$$= H_n(\mathbf{x}(\theta_0^*, \cdot)) - H_n(\mathbf{x}(\hat{\theta}_n^*, \cdot)) + o_p(1)$$
$$\leq H_n(\hat{\mathbf{x}}_n(\theta_0^*, \cdot)) + \sup_{\theta^* \in \Theta_0 \times \Gamma_0} |H_n(\hat{\mathbf{x}}_n(\theta^*, \cdot)) - H_n(\mathbf{x}(\theta^*, \cdot))|$$
$$\qquad - H_n(\hat{\mathbf{x}}_n(\hat{\theta}_n^*, \cdot))$$
$$\qquad + \sup_{\theta^* \in \Theta_0 \times \Gamma_0} |H_n(\hat{\mathbf{x}}_n(\theta^*, \cdot)) - H_n(\mathbf{x}(\theta^*, \cdot))| + o_p(1)$$
$$\leq H_n(\hat{\mathbf{x}}_n(\theta_0^*, \cdot)) - H_n(\hat{\mathbf{x}}_n(\hat{\theta}_n^*, \cdot))$$
$$\qquad + 2 \sup_{\theta^* \in \Theta_0 \times \Gamma_0} |H_n(\hat{\mathbf{x}}_n(\theta^*, \cdot)) - H_n(\mathbf{x}(\theta^*, \cdot))| + o_p(1)$$
$$\leq H_n(\hat{\mathbf{x}}_n(\theta_0^*, \cdot)) - H_n(\hat{\mathbf{x}}_n(\hat{\theta}_n^*, \cdot)) + o_p(1) \qquad \text{by } (5.9).$$

By definition of $\hat{\theta}_n^*$, we have

$$H_n(\hat{\mathbf{x}}_n(\hat{\theta}_n^*, \cdot)) \geq H_n(\hat{\mathbf{x}}_n(\theta_0^*, \cdot))$$

and then $M(\theta_0^*) - M(\hat{\theta}_n^*) \le o_p(1)$. By (5.16) and (5.17), we have

$$[|\hat{\theta}_n^* - \theta_0^*| > \eta] \cap \Pi \subset [M(\theta_0^*) - M(\hat{\theta}_n^*) \ge \gamma] \cap \Pi \subset [o_p(1) \ge \gamma] \cap \Pi.$$

Now

$$P_{\theta_0^*}[|\hat{\theta}_n^* - \theta_0^*| > \eta]$$

$$\le P_{\theta_0^*}([|\hat{\theta}_n^* - \theta_0^*| > \eta] \cap \Pi) + P_{\theta_0^*}(\Pi^c)$$

$$\le P_{\theta_0^*}([o_p(1) \ge \gamma] \cap \Pi) + P_{\theta_0^*}(\Pi^c) \le P_{\theta_0^*}[o_p(1) \ge \gamma] + \varepsilon \qquad \text{by (5.15)}.$$

We have

$$\limsup_{n \to \infty} P_{\theta_0^*}[|\hat{\theta}_n^* - \theta_0^*| > \eta] \le \varepsilon.$$

Since $\varepsilon$ is arbitrary, we have $P_{\theta_0^*}[|\hat{\theta}_n^* - \theta_0^*| > \eta] \to 0$.  $\square$

PROOF OF THEOREM 3.3.   For any convex compact sets $\Theta_0 \subset \Theta$ and $\Gamma_0 \subset \Gamma$ such that $\theta_0^*$ is an interior point of $\Theta_0 \times \Gamma_0$, let $\tilde{\theta}_n^*$ be the maximizer of $H_n(\hat{\mathbf{x}}_n(\theta^*, \cdot))$ in $\Theta_0 \times \Gamma_0$. Then the event

$$(5.18) \qquad\qquad [\tilde{\theta}_n^* \ne \hat{\theta}_n^*] = [\hat{\theta}_n^* \notin \Theta_0 \times \Gamma_0].$$

We will prove the asymptotic normality for $\tilde{\theta}_n^*$ by using Theorem 5.23 in the book of van der Vaart [25]. First, by (5.14), for any $\theta^{*\prime}, \theta^{*\prime\prime} \in \Theta_0 \times \Gamma_0$.

$$(5.19) \qquad |g(y, \mathbf{x}(\theta^{*\prime}, t)) - g(y, \mathbf{x}(\theta^{*\prime\prime}, t))|$$

$$\le \frac{C}{\delta} \left[ 1 + \left| \frac{\partial g}{\partial x}(y, \mathbf{x}(\theta_0^*, t)) \right| \right] |\theta^{*\prime} - \theta^{*\prime\prime}|,$$

where $C$ is a constant depending on $\Theta_0$, $\Gamma_0$ and $T$. $\delta$ is the constant in (5.12). By Assumption 6, $1 + |\frac{\partial g}{\partial x}(y, \mathbf{x}(\theta_0^*, t))|$ has a finite second moment.

Second, we will prove a Taylor expansion for

$$M(\theta^*) = -E_{\theta_0^*}[g(Y_i, \mathbf{x}(\theta^*, T_i))],$$

in the neighborhood of $\theta_0^*$. We expand

$$g(y, \mathbf{x}(\theta^*, t))$$

$$(5.20) \qquad = g(y, \mathbf{x}(\theta_0^*, t))$$

$$+ \left[ \frac{\partial g}{\partial x}(y, \mathbf{x}(\theta_0^*, t)) \right] \left( \frac{\partial \mathbf{x}}{\partial \theta^*}(\theta_0^*, t) \right)^T (\theta^* - \theta_0^*)$$

$$(5.21) \qquad + \frac{1}{2}(\theta^* - \theta_0^*)^T \left[ \frac{\partial g}{\partial x}(y, \mathbf{x}(\theta_0^*, t)) \frac{\partial^2 \mathbf{x}}{\partial \theta^* \partial \theta^{*T}}(\theta_0^*, t) \right.$$

$$\left. + \frac{\partial^2 g}{\partial x^2}(y, \mathbf{x}(\theta_0^*, t)) \frac{\partial \mathbf{x}}{\partial \theta^*}(\theta_0^*, t) \frac{\partial \mathbf{x}}{\partial \theta^*}(\theta_0^*, t)^T \right]$$

$$\times (\theta^* - \theta_0^*) + R_0,$$

where $(\theta^* - \theta_0^*)^T$ denotes the transpose of $(\theta^* - \theta_0^*)$ and $R_0$ is the remainder term. If define a continuous matrix

$$D(y, t, \theta^*) = \frac{\partial g}{\partial x}(y, \mathbf{x}(\theta^*, t)) \frac{\partial^2 \mathbf{x}}{\partial \theta^* \partial \theta^{*T}}(\theta^*, t)$$

$$+ \frac{\partial^2 g}{\partial x^2}(y, \mathbf{x}(\theta^*, t)) \frac{\partial \mathbf{x}}{\partial \theta^*}(\theta^*, t) \frac{\partial \mathbf{x}}{\partial \theta^*}(\theta^*, t)^T,$$

we can express the remainder term by an integral

$$R_0 = (\theta^* - \theta_0^*)^T \left[ \int_0^1 [D(y, t, \theta_0^* + s(\theta^* - \theta_0^*)) - D(y, t, \theta_0^*)](1-s)\,ds \right]$$

$$\times (\theta^* - \theta_0^*).$$

By using the same argument in the proof of Lemma 6, we have

$$|D(y, t, \theta^*)| \leq \frac{C'}{\delta} \left[ 1 + \left| \frac{\partial g}{\partial x}(y, \mathbf{x}(\theta_0^*, t)) \right| + \left| \frac{\partial^2 g}{\partial x^2}(y, \mathbf{x}(\theta_0^*, t)) \right| \right]$$

$$\forall \theta^* \in \Theta_0 \times \Gamma_0,$$

where $C'$ is a constant depending on $\Theta_0$, $\Gamma_0$ and $T$. The right-hand side of the last inequality has a finite expectation by Assumption 6. Hence, by the dominated convergence theorem,

$$E_{\theta_0^*} \left[ \int_0^1 [D(y, t, \theta_0^* + s(\theta^* - \theta_0^*)) - D(y, t, \theta_0^*)](1-s)\,ds \right] \to 0$$

as $\theta^* \to \theta_0^*$. From (5.20), we have

$$(5.22) \qquad M(\theta^*) = M(\theta_0^*) + (\theta^* - \theta_0^*)^T V_{\theta_0^*}(\theta^* - \theta_0^*) + o(\|\theta^* - \theta_0^*\|^2),$$

where

$$V_{\theta_0^*} = -E_{\theta_0^*}[D(y, t, \theta_0^*)]$$

and there is no linear term of $\theta^* - \theta_0^*$ beacuse $\theta_0^*$ is the maximum point of $M$.

Finally, we have

$$H_n(\mathbf{x}(\tilde{\theta}_n^*, \cdot))$$

$$\geq H_n(\hat{\mathbf{x}}_n(\tilde{\theta}_n^*, \cdot)) - \sup_{\theta^* \in \Theta_0 \times \Gamma_0} |H_n(\hat{\mathbf{x}}_n(\theta^*, \cdot)) - H_n(\mathbf{x}_n(\theta^*, \cdot))|$$

$$= \sup_{\theta^* \in \Theta_0 \times \Gamma_0} H_n(\hat{\mathbf{x}}_n(\theta^*, \cdot))$$

$$- \sup_{\theta^* \in \Theta_0 \times \Gamma_0} |H_n(\hat{\mathbf{x}}_n(\theta^*, \cdot)) - H_n(\mathbf{x}_n(\theta^*, \cdot))|$$

$$(5.23)$$

$$\text{by definition of } \tilde{\theta}_n^*$$

$$\geq \sup_{\theta^* \in \Theta_0 \times \Gamma_0} [H_n(\mathbf{x}(\theta^*, \cdot)) - \sup_{\theta^* \in \Theta_0 \times \Gamma_0} |H_n(\mathbf{x}(\theta^*, \cdot)) - H_n(\hat{\mathbf{x}}_n(\theta^*, \cdot))|]$$

$$- \sup_{\theta^* \in \Theta_0 \times \Gamma_0} |H_n(\hat{\mathbf{x}}_n(\theta^*, \cdot)) - H_n(\mathbf{x}_n(\theta^*, \cdot))|$$

$$\geq \sup_{\theta^* \in \Theta_0 \times \Gamma_0} H_n(\mathbf{x}(\theta^*, \cdot)) - 2 \sup_{\theta^* \in \Theta_0 \times \Gamma_0} |H_n(\mathbf{x}(\theta^*, \cdot)) - H_n(\hat{\mathbf{x}}_n(\theta^*, \cdot))|.$$

By (5.9) and

$$\frac{\lambda_n}{n^2} \to \infty \quad \text{and} \quad r_n = o_p\left(\frac{1}{n}\right) \qquad \text{as } n \to \infty,$$

we have

$$\sup_{\theta^* \in \Theta_0 \times \Gamma_0} |H_n(\hat{\mathbf{x}}_n(\theta^*, \cdot)) - H_n(\mathbf{x}(\theta^*, \cdot))|$$

$$= \left[O_p\left(\frac{1}{\sqrt{\lambda_n}}\right)\sqrt{T} + 2T\sqrt{8(8K^2+2)}r_n\right]e^{2KT}O_p(1) + o_p\left(\frac{1}{n}\right)$$

$$= o_p\left(\frac{1}{n}\right).$$

Now if we look $\Theta_0 \times \Gamma_0$ as the parameter space, then by (5.19), (5.22) and (5.23), it follows from Theorem 5.23 in book of van der Vaart [25] that $\sqrt{n}(\hat{\theta}_n^* - \theta_0^*)$ is asymptotically normal with mean zero and covariance matrix

$$(5.24) \quad V_{\theta_0^*}^{-1} E_{\theta_0^*}\left[\left(\frac{\partial g}{\partial x}(y, \mathbf{x}(\theta_0^*, t))\right)^2 \left(\frac{\partial \mathbf{x}}{\partial \theta^*}(\theta_0^*, t)\right)^T \left(\frac{\partial \mathbf{x}}{\partial \theta^*}(\theta_0^*, t)\right)\right] V_{\theta_0^*}^{-1}.$$

Note the asymptotic covariance matrix does not depend on $\Theta_0 \times \Gamma_0$. Because $\hat{\theta}_n^*$'s are tight and by (5.18), we can make

$$\sup_n P_{\theta_0^*}[\tilde{\theta}_n^* \neq \hat{\theta}_n^*] = \sup_n P_{\theta_0^*}[\hat{\theta}_n^* \notin \Theta_0 \times \Gamma_0]$$

arbitrarily small by taking large $\Theta_0 \times \Gamma_0$. It follows the lemma below that $\sqrt{n}(\hat{\theta}_n^* - \theta_0^*)$ is asymptotically normal with mean zero and covariance matrix (5.24). Similarly, we can prove the asymptotic normality of $\hat{\tilde{\theta}}_n^*$ with the same asymptotic covariance matrix. $\square$

LEMMA 8. *Let $\{X_n : n = 1, 2, \ldots\}$ be a sequence of random variables. For each $m = 1, 2, \ldots$, there is a sequence of random variables $\{X_n^{(m)} : n = 1, 2, \ldots\}$ such that for any $\varepsilon > 0$, we have*

$$\lim_{m \to \infty} \sup_n P(|X_n - X_n^{(m)}| > \varepsilon) = 0.$$

*Suppose that for each $m$, the sequence $\{X_n^{(m)} : n = 1, 2, \ldots\}$ converges weakly to the same random variable $X$ which does not depend on $m$. Then $\{X_n : n = 1, 2, \ldots\}$ converges weakly to $X$.*

Proof.  We calculate the characteristic function of $\{X_n : n = 1, 2, \ldots\}$. Fix $u \in \mathbb{R}$.

$$|E(e^{iuX_n}) - E(e^{iuX_n^{(m)}})| \leq E|e^{iuX_n} - e^{iuX_n^{(m)}}| = E|e^{iu(X_n - X_n^{(m)})} - 1|.$$

Because $e^{iut}$ is a continuous function of $t$ and bounded by 2, for any $\delta > 0$, there exists $\varepsilon > 0$ such that

$$|e^{iut} - 1| \leq \delta \qquad \forall |t| \leq \varepsilon.$$

We have

$$|E(e^{iuX_n}) - E(e^{iuX_n^{(m)}})|$$
$$\leq E|e^{iu(X_n - X_n^{(m)})} - 1|$$
$$= E[|e^{iu(X_n - X_n^{(m)})} - 1|\mathbf{1}_{[|X_n - X_n^{(m)}| > \varepsilon]}] + E[|e^{iu(X_n - X_n^{(m)})} - 1|\mathbf{1}_{[|X_n - X_n^{(m)}| \leq \varepsilon]}]$$
$$\leq 2P(|X_n - X_n^{(m)}| > \varepsilon) + \delta.$$

Hence,

$$\lim_{m \to \infty} \sup_n |E(e^{iuX_n}) - E(e^{iuX_n^{(m)}})| \leq \delta.$$

Since $\delta$ is arbitrary, we have

$$\lim_{m \to \infty} \sup_n |E(e^{iuX_n}) - E(e^{iuX_n^{(m)}})| = 0.$$

Because for each $m$, the sequence $\{X_n^{(m)} : n = 1, 2, \ldots\}$ converges weakly to $X$, $E(e^{iuX_n^{(m)}}) \to E(e^{iuX})$. So we have $E(e^{iuX_n}) \to E(e^{iuX})$, that is, $\{X_n : n = 1, 2, \ldots\}$ converges weakly to $X$.  □

Proof of Lemma 5.  If the lemma is wrong, then there exist $M > 0$, $\delta > 0$ and a sequence $\{\mathbb{L}_q : q \geq 1\}$ with $r_q \to 0$, such that for each $q$, there exist

$$(\mathbf{w}_q, \mathbf{v}_q) \in \operatorname*{arg\,min}_{\substack{\hat{\mathbf{x}}, \hat{\mathbf{z}} \in \mathbb{L}_q, \hat{\mathbf{x}}(0) = x, \hat{\mathbf{z}} = z \\ \|\dot{\hat{\mathbf{x}}}\|_\infty \leq M, \|\dot{\hat{\mathbf{z}}}\|_\infty \leq M}} J(\hat{\mathbf{x}}, \hat{\mathbf{z}}, \theta)$$

with

(5.25) $$\|\mathbf{w}_q - \mathbf{x}\|_\infty \geq \delta \quad \text{or} \quad \|\mathbf{v}_q - \mathbf{z}\|_\infty \geq \delta.$$

We will show that $\{(\mathbf{w}_q, \mathbf{v}_q) : q \geq 1\}$ are equicontinuous. Fix any $\eta > 0$. Let $K$ be a positive constant such that

$$|F(x, z, t, \theta)| \leq K, \qquad |G(x, z, t, \theta)| \leq K \qquad \forall |x|, |z| \leq M, t \in [0, T]$$

For any $t_0 \in [0, T]$,

$$
\begin{aligned}
|\mathbf{w}_q(t) - \mathbf{w}_q(t_0)| &= \left| \int_{t_0}^{t} \frac{d\mathbf{w}_q}{dt}(s) \, ds \right| \leq \int_{t_0}^{t} \left| \frac{d\mathbf{w}_q}{dt}(s) \right| ds \\
&= \int_{t_0}^{t} \left| \frac{d\mathbf{w}_q}{dt}(s) - F(\mathbf{w}_q, \mathbf{v}_q, t, \theta) + F(\mathbf{w}_q, \mathbf{v}_q, t, \theta) \right| ds \\
&\leq \int_{t_0}^{t} \left| \frac{d\mathbf{w}_q}{dt}(s) - F(\mathbf{w}_q, \mathbf{v}_q, t, \theta) \right| ds + K|t - t_0| \\
&\leq \sqrt{\int_{t_0}^{t} \left| \frac{d\mathbf{w}_q}{dt}(s) - F(\mathbf{w}_q, \mathbf{v}_q, t, \theta) \right|^2 ds} + K|t - t_0| \\
&\leq \sqrt{J(\mathbf{w}_q, \mathbf{v}_q, \theta)} + K|t - t_0|.
\end{aligned}
$$

By the similar argument as in the proof of Theorem 3.1, we have

$$
J(\mathbf{w}_q, \mathbf{v}_q, \theta) \to 0,
$$

because $r_q \to 0$. Therefore, there exists $q_0$ such that if $q \geq q_0$, $J(\mathbf{w}_q, \mathbf{v}_q, \theta) \leq (\frac{\eta}{2})^2$. So for any $q \geq q_0$ and $|t - t_0| \leq \frac{\eta}{2K}$, we have

$$
|\mathbf{w}_q(t) - \mathbf{w}_q(t_0)| \leq \eta.
$$

Hence, $\{\mathbf{w}_q : q \geq 1\}$ are equicontinuous. Similarly, $\{\mathbf{v}_q : q \geq 1\}$ are equicontinuous. Then by Ascoli's theorem there is a uniformly convergent subsequence. Without loss of generality, we assume that $\mathbf{w}_q \to \mathbf{x}_0$ and $\mathbf{v}_q \to \mathbf{z}_0$ uniformly. We will show that $(\mathbf{x}_0, \mathbf{z}_0)$ are also the solutions of (1.1).

$$
\int_{t_0}^{t} \left| \frac{d\mathbf{w}_q}{dt}(s) - F(\mathbf{w}_q, \mathbf{v}_q, t, \theta) \right| ds \leq \sqrt{J(\mathbf{w}_q, \mathbf{v}_q, \theta)} \to 0.
$$

Hence, $\frac{d\mathbf{w}_q}{dt} \to F(\mathbf{x}_0, \mathbf{z}_0, t, \theta)$ in $L^1[0, T]$. Then we have

$$
\mathbf{w}_q \to x + \int_{0}^{t} F(\mathbf{x}_0(s), \mathbf{z}_0(s), s, \theta) \, ds
$$

uniformly, that is,

$$
\mathbf{x}_0(t) = x + \int_{0}^{t} F(\mathbf{x}_0(s), \mathbf{z}_0(s), s, \theta) \, ds.
$$

Similarly

$$
\mathbf{z}_0(t) = z + \int_{0}^{t} G(\mathbf{x}_0(s), \mathbf{z}_0(s), s, \theta) \, ds.
$$

$(\mathbf{x}_0, \mathbf{z}_0)$ are the solutions of (1.1). By the uniqueness, $(\mathbf{x}, \mathbf{z}) = (\mathbf{x}_0, \mathbf{z}_0)$, so $(\mathbf{w}_q, \mathbf{v}_q) \to (\mathbf{x}_0, \mathbf{z}_0)$ uniformly. This contradicts (5.25). $\quad\square$

PROOF OF THEOREM 4.1.   Let $\tau = (0 = t_1 < \cdots < t_{l+1} = T)$. Then the dimension of $\mathbb{L}$ is $l+3$. Let $\{\phi_j : j = -2, -1, 0, \ldots l\}$ be the B-spline bases. For each $1 \leq j \leq l-3$, $\phi_j$ is a piecewise polynomial of order 4 and vanishes outside the interval $(t_j, t_{j+4})$. All the bases are bounded by 1. For any $\hat{\mathbf{x}} \in \mathbb{L}$, let

$$\hat{\mathbf{x}} = \sum_{j=-2}^{l+1} c_j \phi_j,$$

where $(c_j, j = -2, \ldots, l+1)$ are coefficients.

Because $F$ has the third-order continuous partial derivatives, $\mathbf{x}$ has a continuous fourth derivative. Let

$$R = \|\mathbf{x}\|_\infty, \qquad R_1 = \left\|\frac{d\mathbf{x}}{dt}\right\|_\infty, \qquad K_2 = \left\|\frac{d^2\mathbf{x}}{dt^2}\right\|_\infty \quad \text{and} \quad K_4 = \left\|\frac{d^4\mathbf{x}}{dt^4}\right\|_\infty.$$

We can also find positive number $K_0$, $K$ and $K_1$, such that

$$(5.26) \quad \begin{aligned} |F(x,z,t)| &\leq K_0, & |F(x,t) - F(x',t)| &\leq K|x - x'|, \\ |F_x(x,t)| &\leq K_1, & |F_t(x,t)| &\leq K_1 & \forall |x| \leq R+1, |x'| \leq R+1. \end{aligned}$$

By definition of $r$ and the proof of Lemma 1, we have

$$(5.27) \quad r \leq \max\{C_0 K_4 |\tau|^4, C_1 K_4 |\tau|^3\} = C_1 K_4 |\tau|^3,$$

where $C_0 = \frac{5}{384} < C_1 = \frac{9+\sqrt{3}}{216}$. Because $F_x(\mathbf{x}(t), t) < 0$ for all $0 \leq t \leq T$, we can find positive numbers $\delta_2 < 1$ and $\gamma$ such that

$$(5.28) \quad F_x(x,t) \leq -\gamma \qquad \forall |x - \mathbf{x}(t)| \leq \delta_2, 0 \leq t \leq T.$$

Let

$$\hat{\mathbf{x}}_0 \in \underset{\hat{\mathbf{x}} \in \mathbb{L}, \hat{\mathbf{x}}(0) = x, \|\hat{\mathbf{x}} - \mathbf{x}\|_\infty < \delta_2}{\arg\min} J(\hat{\mathbf{x}}).$$

First, we show:

LEMMA 9.   If $2r < \delta_2$, we have

$$J(\hat{\mathbf{x}}_0) \leq 4(4K^2 + 2)\operatorname{Tr}^2.$$

PROOF.   By the definition (4.5) of $r$, there exist $\mathbf{w} \in \mathbb{L}$ with $\mathbf{w}(0) = x$ such that

$$\|\mathbf{w} - \mathbf{x}\|_\infty \leq 2r < \delta_2 < 1, \qquad \left\|\frac{d\mathbf{w}}{dt} - \frac{d\mathbf{x}}{dt}\right\|_\infty \leq 2r.$$

By (5.26), we have

$$\int_0^T \left| \frac{d\mathbf{w}}{dt} - F(\mathbf{w}, t) \right|^2 dt$$

$$= \int_0^T \left| \frac{d\mathbf{w}}{dt} - \frac{d\mathbf{x}}{dt} + F(\mathbf{x}, t) - F(\mathbf{w}, t) \right|^2 dt$$

$$\leq 2 \int_0^T \left| \frac{d\mathbf{w}}{dt} - \frac{d\mathbf{x}}{dt} \right|^2 dt + 2 \int_0^T |F(\mathbf{x}, t) - F(\mathbf{w}, t)|^2 dt$$

$$\leq 2 \int_0^T \left| \frac{d\mathbf{w}}{dt} - \frac{d\mathbf{x}}{dt} \right|^2 dt + 4K^2 \int_0^T |\mathbf{x} - \mathbf{w}|^2 dt \leq 4(4K^2 + 2) \operatorname{Tr}^2.$$

By the definition of $\hat{\mathbf{x}}_0$, we have

$$J(\hat{\mathbf{x}}_0) \leq J(\mathbf{w}) = \int_0^T \left| \frac{d\mathbf{w}}{dt} - F(\mathbf{w}, t) \right|^2 dt \leq 4(4K^2 + 2) \operatorname{Tr}^2. \qquad \square$$

LEMMA 10. *For any $\bar{t} \in \{0 < s < T : \frac{d\hat{\mathbf{x}}_0}{dt}(s) - \frac{d\mathbf{x}}{dt}(s) = 0\}$,*

$$|\mathbf{x}(\bar{t}) - \hat{\mathbf{x}}_0(\bar{t})| \leq \beta_1 \kappa \left\| \frac{d^2 \hat{\mathbf{x}}_0}{dt^2} \right\|_{L^2[0,T]} |\tau|^{1/2} + \beta_2 \kappa |\tau|$$

$$+ \kappa (4\sqrt{6}\kappa + \beta_3) \beta_4 \sqrt{T} |\tau|^{3/2} + \beta_6 \sqrt{T} |\tau|^{7/2},$$

*where $\beta_1, \beta_2, \beta_3, \beta_4, \beta_6$ are constants depending only on $\mathbf{x}$ and $F$.*

PROOF. Pick $i$ such that $\bar{t} \in (t_i, t_{i+4}) \subset [0, T]$. There exists a small positive number $\eta$, such that for any $|u| \leq \eta$, $u \in \mathbb{R}$, we have

$$\|(\hat{\mathbf{x}}_0 + u\phi_i) - \mathbf{x}\|_\infty < \delta_2.$$

By the definition of $\hat{\mathbf{x}}_0$, we have

$$J(\hat{\mathbf{x}}_0) \leq J(\hat{\mathbf{x}}_0 + u\phi_i),$$

that is, 0 is the minimum point of the function $J(\hat{\mathbf{x}}_0 + u\phi_i)$ of $u$ in the region $\{u : |u| \leq \eta, u \in \mathbb{R}\}$. Therefore,

$$0 = \left. \frac{\partial J}{\partial u}(\hat{\mathbf{x}}_0 + u\phi_i) \right|_{u=0}$$

$$= \int_{t_i}^{t_{i+4}} \left( \frac{d\phi_i}{dt} - F_x(\hat{\mathbf{x}}_0(t), t)\phi_i(t) \right) \left( \frac{d\hat{\mathbf{x}}_0}{dt} - F(\hat{\mathbf{x}}_0(t), t) \right) dt$$

$$= \int_{t_i}^{t_{i+4}} \frac{d\phi_i}{dt} \left( \frac{d\hat{\mathbf{x}}_0}{dt} - F(\hat{\mathbf{x}}_0(t), t) \right) dt$$

$$-\int_{t_i}^{t_{i+4}} F_x(\hat{\mathbf{x}}_0(t),t)\phi_i(t)\left(\frac{d\hat{\mathbf{x}}_0}{dt} - \frac{d\mathbf{x}}{dt} + F(\mathbf{x}(t),t) - F(\hat{\mathbf{x}}_0(t),t)\right)dt$$

$$(5.29) \quad = \int_{t_i}^{t_{i+4}} \frac{d\phi_i}{dt}\left(\frac{d\hat{\mathbf{x}}_0}{dt} - F(\hat{\mathbf{x}}_0(t),t)\right)dt$$

$$-\int_{t_i}^{t_{i+4}} F_x(\hat{\mathbf{x}}_0(t),t)\phi_i(t)\left(\frac{d\hat{\mathbf{x}}_0}{dt} - \frac{d\mathbf{x}}{dt}\right)$$

$$-\int_{t_i}^{t_{i+4}} F_x(\hat{\mathbf{x}}_0(t),t)\phi_i(t)(F(\mathbf{x}(t),t) - F(\mathbf{x}(\bar{t}),\bar{t})$$

$$+ F(\hat{\mathbf{x}}_0(\bar{t}),\bar{t}) - F(\hat{\mathbf{x}}_0(t),t))\,dt$$

$$-\int_{t_i}^{t_{i+4}} F_x(\hat{\mathbf{x}}_0(t),t)\phi_i(t)(F(\mathbf{x}(\bar{t}),\bar{t}) - F(\hat{\mathbf{x}}_0(\bar{t}),\bar{t}))\,dt.$$

The integrals are from $t_j$ to $t_{j+4}$ because $\phi_i$ and $\frac{d\phi_i}{dt}$ vanish outside $(t_j, t_{j+4})$. We will estimate every term on the right-hand sides of (5.29).

First, by the formula (see DeBoor [8])

$$\frac{d\phi_i}{dt} = \frac{3}{t_{i+3} - t_i}\phi_{i,3} - \frac{3}{t_{i+4} - t_{i+1}}\phi_{i+1,3},$$

where $\phi_{i,3}$'s are the B-spline bases of order 3 with knots $\tau$, we can calculate

$$\left|\int_{t_i}^{t_{i+4}} \frac{d\phi_i}{dt}\left(\frac{d\hat{\mathbf{x}}_0}{dt} - F(\hat{\mathbf{x}}_0,\hat{\mathbf{z}}_0,t)\right)dt\right|$$

$$(5.30)$$

$$\leq \left\|\frac{d\phi_i}{dt}\right\|_{L^2[0,T]} \left\|\frac{d\hat{\mathbf{x}}_0}{dt} - F(\hat{\mathbf{x}}_0,\hat{\mathbf{z}}_0,t)\right\|_{L^2[0,T]}$$

$$\leq \left(\frac{3}{\sqrt{t_{i+3} - t_i}} + \frac{3}{\sqrt{t_{i+4} - t_{i+1}}}\right)\sqrt{J(\hat{\mathbf{x}}_0,\hat{\mathbf{z}}_0,t)}$$

$$\leq \frac{6}{\sqrt{3\min_j |t_{j+1} - t_j|}}\sqrt{J(\hat{\mathbf{x}}_0,\hat{\mathbf{z}}_0,t)}$$

$$= \frac{6}{\sqrt{3|\tau|/\kappa}}\sqrt{J(\hat{\mathbf{x}}_0,\hat{\mathbf{z}}_0,t)} \leq \frac{12}{\sqrt{3|\tau|/\kappa}}\sqrt{(4K^2+2)T}C_1 K_4 |\tau|^3$$

$$\leq 4\sqrt{6(4K^2+2)T\kappa}\,C_1 K_4 |\tau|^{5/2}.$$

Because $\frac{d\hat{\mathbf{x}}_0}{dt}(\bar{t}) - \frac{d\mathbf{x}}{dt}(\bar{t}) = 0$,

$$\left|\int_{t_i}^{t_{i+4}} F_x(\hat{\mathbf{x}}_0(t),t)\phi_i(t)\left(\frac{d\hat{\mathbf{x}}_0}{dt}(t) - \frac{d\mathbf{x}}{dt}(t)\right)dt\right|$$

$$= \left|\int_{t_i}^{t_{i+4}} F_x(\hat{\mathbf{x}}_0(t),t)\phi_i(t)\left(\frac{d\hat{\mathbf{x}}_0}{dt}(t) - \frac{d\mathbf{x}}{dt}(t) - \frac{d\hat{\mathbf{x}}_0}{dt}(\bar{t}) + \frac{d\mathbf{x}}{dt}(\bar{t})\right)dt\right|$$

$$\leq K_1 \int_{t_i}^{t_{i+4}} \left| \frac{d\hat{\mathbf{x}}_0}{dt}(t) - \frac{d\hat{\mathbf{x}}_0}{dt}(\bar{t}) \right| dt + K_1 \int_{t_i}^{t_{i+4}} \left| \frac{d\mathbf{x}}{dt}(t) - \frac{d\mathbf{x}}{dt}(\bar{t}) \right| dt$$

$$\leq K_1 \int_{t_i}^{t_{i+4}} \left| \int_{\bar{t}}^{t} \frac{d^2\hat{\mathbf{x}}_0}{dt^2}(s)\, ds \right| dt + K_1 \int_{t_i}^{t_{i+4}} \left| \int_{\bar{t}}^{t} \frac{d^2\mathbf{x}}{dt^2}(s)\, ds \right| dt$$

$$\leq K_1 \int_{t_i}^{t_{i+4}} \int_{t_i}^{t} \left| \frac{d^2\hat{\mathbf{x}}_0}{dt^2}(s) \right| ds\, dt + K_1 \int_{t_i}^{t_{i+4}} \int_{t_i}^{t} \left| \frac{d^2\mathbf{x}}{dt^2}(s) \right| ds\, dt$$

$$(5.31) \quad \leq K_1 \int_{t_i}^{t_{i+4}} \left| \frac{d^2\hat{\mathbf{x}}_0}{dt^2}(s) \right| \int_{s}^{t_{i+4}} dt\, ds + K_1 \int_{t_i}^{t_{i+4}} \int_{t_i}^{t} K_2\, ds\, dt$$

$$= K_1 \int_{t_i}^{t_{i+4}} \left| \frac{d^2\hat{\mathbf{x}}_0}{dt^2}(s) \right| (t_{i+4} - s)\, ds + \frac{K_1 K_2}{2}(t_{i+4} - t_i)^2$$

$$\leq K_1 \left\| \frac{d^2\hat{\mathbf{x}}_0}{dt^2} \right\|_{L^2[0,T]} \left[ \int_{t_i}^{t_{i+4}} (t_{i+4} - s)^2\, ds \right]^{1/2}$$
$$\quad + \frac{K_1 K_2}{2}(t_{i+4} - t_i)^2$$

$$\leq \frac{K_1}{\sqrt{3}} \left\| \frac{d^2\hat{\mathbf{x}}_0}{dt^2} \right\|_{L^2[0,T]} (t_{i+4} - t_i)^{3/2} + \frac{K_1 K_2}{2}(t_{i+4} - t_i)^2$$

$$\leq \frac{8K_1}{\sqrt{3}} \left\| \frac{d^2\hat{\mathbf{x}}_0}{dt^2} \right\|_{L^2[0,T]} |\tau|^{3/2} + 8K_1 K_2 |\tau|^2.$$

For the second term from the last on the right-hand side of (5.29),

$$\left| \int_{t_i}^{t_{i+4}} F_x(\hat{\mathbf{x}}_0(t), t)\phi_i(t)(F(\mathbf{x}(t), t) - F(\mathbf{x}(\bar{t}), \bar{t}))\, dt \right|$$

$$\leq K_1 \int_{t_i}^{t_{i+4}} |F(\mathbf{x}(t), t) - F(\mathbf{x}(\bar{t}), \bar{t})|\, dt$$

$$(5.32) \quad = K_1 \int_{t_i}^{t_{i+4}} \left| \int_{\bar{t}}^{t} \left( F_x(\mathbf{x}(s), s)\frac{d\mathbf{x}}{dt} + F_t(\mathbf{x}(s), s) \right) ds \right| dt$$

$$\leq K_1 \int_{t_i}^{t_{i+4}} \int_{t_i}^{t} (K_1 R_1 + K_1)\, ds\, dt \leq K_1 \frac{K_1 R_1 + K_1}{2}(t_{i+4} - t_i)^2$$

$$\leq 8K_1^2(R_1 + 1)|\tau|^2$$

and

$$\left| \int_{t_i}^{t_{i+4}} F_x(\hat{\mathbf{x}}_0(t), t)\phi_i(t)(F(\hat{\mathbf{x}}_0(\bar{t}), \bar{t}) - F(\hat{\mathbf{x}}_0(t), t))\, dt \right|$$

$$\leq K_1 \int_{t_i}^{t_{i+4}} |F(\hat{\mathbf{x}}_0(t), t) - F(\hat{\mathbf{x}}_0(\bar{t}), \bar{t})| \, dt$$

$$= K_1 \int_{t_i}^{t_{i+4}} \left| \int_{\bar{t}}^{t} \left( F_x(\hat{\mathbf{x}}_0(s), s) \frac{d\hat{\mathbf{x}}_0}{dt} + F_t(\hat{\mathbf{x}}_0(s), s) \right) ds \right| dt$$

$$\leq K_1 \int_{t_i}^{t_{i+4}} \int_{t_i}^{t} \left( K_1 \left| \frac{d\hat{\mathbf{x}}_0}{dt}(s) \right| + K_1 \right) ds \, dt$$

$$\leq K_1^2 \int_{t_i}^{t_{i+4}} \int_{t_i}^{t} \left| \frac{d\hat{\mathbf{x}}_0}{dt}(s) \right| ds \, dt + \frac{K_1^2}{2}(t_{i+4} - t_i)^2$$

$$= K_1^2 \int_{t_i}^{t_{i+4}} \left| \frac{d\hat{\mathbf{x}}_0}{dt}(s) \right| (t_{i+4} - s) \, ds \, dt + \frac{K_1^2}{2}(t_{i+4} - t_i)^2$$

$$= K_1^2 \int_{t_i}^{t_{i+4}} \left| \frac{d\hat{\mathbf{x}}_0}{dt}(s) - F(\hat{\mathbf{x}}_0(s), s) + F(\hat{\mathbf{x}}_0(s), s) \right| (t_{i+4} - s) \, ds \, dt$$

$$+ \frac{K_1^2}{2}(t_{i+4} - t_i)^2$$

$$= K_1^2 \int_{t_i}^{t_{i+4}} \left| \frac{d\hat{\mathbf{x}}_0}{dt}(s) - F(\hat{\mathbf{x}}_0(s), s) \right| (t_{i+4} - s) \, ds \, dt$$

$$+ \frac{K_1^2 K_0}{2}(t_{i+4} - t_i)^2 + \frac{K_1^2}{2}(t_{i+4} - t_i)^2$$

(5.33)
$$\leq \frac{K_1^2}{\sqrt{3}} \left\| \frac{d^2 \hat{\mathbf{x}}_0}{dt^2} - F(\hat{\mathbf{x}}_0, \cdot) \right\|_{L^2[0,T]} (t_{i+4} - t_i)^{3/2}$$

$$+ \frac{K_1^2}{2}(K_0 + 1)(t_{i+4} - t_i)^2$$

$$\leq \frac{8K_1^2}{\sqrt{3}} \left\| \frac{d^2 \hat{\mathbf{x}}_0}{dt^2} - F(\hat{\mathbf{x}}_0, \cdot) \right\|_{L^2[0,T]} |\tau|^{3/2} + 8K_1^2(K_0 + 1)|\tau|^2$$

$$= \frac{8K_1^2}{\sqrt{3}} \sqrt{J(\hat{\mathbf{x}}_0, \hat{\mathbf{z}}_0, t)} |\tau|^{3/2} + 8K_1^2(K_0 + 1)|\tau|^2$$

$$\leq \frac{16K_1^2}{\sqrt{3}} \sqrt{(4K^2 + 2)T} C_1 K_4 |\tau|^{9/2} + 8K_1^2(K_0 + 1)|\tau|^2.$$

Now we calculate the last term,

$$\left| \int_{t_i}^{t_{i+4}} F_x(\hat{\mathbf{x}}_0(t), t) \phi_i(t)(F(\mathbf{x}(\bar{t}), \bar{t}) - F(\hat{\mathbf{x}}_0(\bar{t}), \bar{t})) \, dt \right|$$

$$= \left| \int_{t_i}^{t_{i+4}} F_x(\hat{\mathbf{x}}_0(t), t) \phi_i(t) \, dt \right| |F(\mathbf{x}(\bar{t}), \bar{t}) - F(\hat{\mathbf{x}}_0(\bar{t}), \bar{t})|$$

$$(5.34) \qquad = \left| \int_{t_i}^{t_{i+4}} F_x(\hat{\mathbf{x}}_0(t), t) \phi_i(t) \, dt \right| |F_x(x', \bar{t})(\mathbf{x}(\bar{t}) - \hat{\mathbf{x}}_0(\bar{t}))|$$

$$\geq \gamma \left| \int_{t_i}^{t_{i+4}} \phi_i(t) \, dt \right| \cdot \gamma |\mathbf{x}(\bar{t}) - \hat{\mathbf{x}}_0(\bar{t})|$$

$$= \frac{\gamma^2 |t_{i+4} - t_i|}{4} |\mathbf{x}(\bar{t}) - \hat{\mathbf{x}}_0(\bar{t})|$$

$$\geq \frac{\gamma^2 |\tau|}{\kappa} |\mathbf{x}(\bar{t}) - \hat{\mathbf{x}}_0(\bar{t})|,$$

where $x'$ is a number between $\mathbf{x}(\bar{t})$ and $\hat{\mathbf{x}}_0(\bar{t})$ and we use the formula (see equality (4.29) in Schumaker [23] and Theorem 4.23 in Schumaker [23])

$$\int_{t_i}^{t_{i+4}} \phi_i(t) \, dt = \frac{|t_{i+4} - t_i|}{4}.$$

From (5.29)–(5.34), we have

$$\frac{\gamma^2 |\tau|}{\kappa} |\mathbf{x}(\bar{t}) - \hat{\mathbf{x}}_0(\bar{t})|$$

$$\leq \left| \int_{t_i}^{t_{i+4}} F_x(\hat{\mathbf{x}}_0(t), t) \phi_i(t) (F(\mathbf{x}(\bar{t}), \bar{t}) - F(\hat{\mathbf{x}}_0(\bar{t}), \bar{t})) \, dt \right|$$

$$\leq \left| \int_{t_i}^{t_{i+4}} \frac{d\phi_i}{dt} \left( \frac{d\hat{\mathbf{x}}_0}{dt} - F(\hat{\mathbf{x}}_0(t), t) \right) dt \right|$$

$$\quad + \left| \int_{t_i}^{t_{i+4}} F_x(\hat{\mathbf{x}}_0(t), t) \phi_i(t) \left( \frac{d\hat{\mathbf{x}}_0}{dt} - \frac{d\mathbf{x}}{dt} \right) \right|$$

$$\quad + \left| \int_{t_i}^{t_{i+4}} F_x(\hat{\mathbf{x}}_0(t), t) \phi_i(t) (F(\mathbf{x}(t), t) - F(\mathbf{x}(\bar{t}), \bar{t})) \, dt \right|$$

$$\quad + \left| \int_{t_i}^{t_{i+4}} F_x(\hat{\mathbf{x}}_0(t), t) \phi_i(t) (F(\hat{\mathbf{x}}_0(\bar{t}), \bar{t}) - F(\hat{\mathbf{x}}_0(t), t)) \, dt \right|$$

$$\leq 4\sqrt{6(4K^2 + 2)T} \kappa C_1 K_4 |\tau|^{5/2}$$

$$\quad + \left( \frac{8K_1}{\sqrt{3}} \left\| \frac{d^2\hat{\mathbf{x}}_0}{dt^2} \right\|_{L^2[0,T]} |\tau|^{3/2} + 8K_1 K_2 |\tau|^2 \right)$$

$$\quad + 8K_1^2 (R_1 + 1) |\tau|^2$$

$$\quad + \left( \frac{8K_1^2}{\sqrt{3}} \sqrt{(4K^2 + 2)T} C_1 K_4 |\tau|^{9/2} + 8K_1^2 (K_0 + 1) |\tau|^2 \right)$$

$$\leq \frac{8K_1}{\sqrt{3}} \left\| \frac{d^2\hat{\mathbf{x}}_0}{dt^2} \right\|_{L^2[0,T]} |\tau|^{3/2}$$

$$+ (8K_1K_2 + 8K_1^2(R_1 + 1) + 8K_1^2(K_0 + 1))|\tau|^2$$
$$+ \left(4\sqrt{6\kappa} + \frac{16K_1^2}{\sqrt{3}}\right)\sqrt{(4K^2 + 2)T}C_1K_4|\tau|^{5/2}$$
$$+ \frac{8K_1^2}{\sqrt{3}}\sqrt{(4K^2 + 2)T}C_1K_4|\tau|^{9/2}.$$

Therefore,

$$|\mathbf{x}(\bar{t}) - \hat{\mathbf{x}}_0(\bar{t})| \leq \beta_1\kappa\left\|\frac{d^2\hat{\mathbf{x}}_0}{dt^2}\right\|_{L^2[0,T]}|\tau|^{1/2} + \beta_2\kappa|\tau|$$
$$+ \kappa(4\sqrt{6\kappa} + \beta_3)\beta_4\sqrt{T}|\tau|^{3/2} + \beta_6\sqrt{T}|\tau|^{7/2},$$

where

$$\beta_1 = \frac{8K_1}{\sqrt{3}\gamma^2},$$
$$\beta_2 = \frac{1}{\gamma^2}(8K_1K_2 + 8K_1^2(R_1 + 1) + 8K_1^2(K_0 + 1)),$$
$$\beta_3 = \frac{16K_1^2}{\sqrt{3}}, \qquad \beta_4 = \frac{1}{\gamma^2}\sqrt{(4K^2 + 2)}C_1K_4,$$
$$\beta_6 = \frac{8K_1^2}{\sqrt{3}}\sqrt{(4K^2 + 2)}C_1K_4. \qquad \square$$

Let

$$t_0 \in \arg\max_{0 \leq s \leq T} |\mathbf{x}(s) - \hat{\mathbf{x}}_0(s)|.$$

Then either $t_0 \in (0, T)$ or $t_0 = T$. If $t_0 \in (0, T)$, we have $\frac{d\hat{\mathbf{x}}_0}{dt}(t_0) - \frac{d\mathbf{x}}{dt}(t_0) = 0$. By Lemma 10,

$$(5.35) \quad \begin{aligned} \|\mathbf{x} - \hat{\mathbf{x}}_0\|_\infty &= |\mathbf{x}(t_0) - \hat{\mathbf{x}}_0(t_0)| \\ &\leq \beta_1\kappa\left\|\frac{d^2\hat{\mathbf{x}}_0}{dt^2}\right\|_{L^2[0,T]}|\tau|^{1/2} + \beta_2\kappa|\tau| \\ &\quad + \kappa(4\sqrt{6\kappa} + \beta_3)\beta_4\sqrt{T}|\tau|^{3/2} + \beta_6\sqrt{T}|\tau|^{7/2}. \end{aligned}$$

In the case $t_0 = T$, let

$$s_0 = \inf\left\{0 \leq s \leq T : \frac{d\hat{\mathbf{x}}_0}{dt}(t) - \frac{d\mathbf{x}}{dt}(t) \neq 0, \forall s < t < T\right\} \wedge T.$$

If $s_0 = T$, there is an increasing sequence $\{s_n, n \geq 1\}$ which converges to $T$ and $\frac{d\hat{\mathbf{x}}_0}{dt}(s_n) - \frac{d\mathbf{x}}{dt}(s_n) = 0$ for all $n$, hence we still have inequality (5.35). If $s_0 < T$, then either $s_0 = 0$ or $\frac{d\hat{\mathbf{x}}_0}{dt}(s_0) - \frac{d\mathbf{x}}{dt}(s_0) = 0$. In both cases, we have

$$|\mathbf{x}(s_0) - \hat{\mathbf{x}}_0(s_0)|$$
$$\leq \beta_1 \kappa \left\| \frac{d^2\hat{\mathbf{x}}_0}{dt^2} \right\|_{L^2[0,T]} |\tau|^{1/2} + \beta_2 \kappa |\tau|$$
$$+ \kappa(4\sqrt{6\kappa} + \beta_3)\beta_4\sqrt{T}|\tau|^{3/2} + \beta_6\sqrt{T}|\tau|^{7/2}.$$

Without loss of generality, we assume that $\frac{d\hat{\mathbf{x}}_0}{dt}(t) - \frac{d\mathbf{x}}{dt}(t) > 0$ for all $t \in (s_0, T)$. Hence, $\hat{\mathbf{x}}_0(t) - \mathbf{x}(t)$ is increasing in $(s_0, T)$ and $\hat{\mathbf{x}}_0(T) - \mathbf{x}(T) > 0$. We have the following two cases,

- If $\hat{\mathbf{x}}_0(s_0) - \mathbf{x}(s_0) \geq 0$, then for any $t \in (s_0, T)$,

$$F(\mathbf{x}(t), t) - F(\hat{\mathbf{x}}_0(t), t) = F_x(x', t)(\mathbf{x}(t) - \hat{\mathbf{x}}_0(t)) \geq 0.$$

Now we have

$$J(\hat{\mathbf{x}}_0) \geq \int_{s_0}^{T} \left( \frac{d\hat{\mathbf{x}}_0}{dt} - F(\hat{\mathbf{x}}_0(t), t) \right)^2 dt$$
$$= \int_{s_0}^{T} \left( \frac{d\hat{\mathbf{x}}_0}{dt} - \frac{d\mathbf{x}}{dt} + F(\mathbf{x}(t), t) - F(\hat{\mathbf{x}}_0(t), t) \right)^2 dt$$
$$\geq \int_{s_0}^{T} \left( \frac{d\hat{\mathbf{x}}_0}{dt} - \frac{d\mathbf{x}}{dt} \right)^2 dt$$
$$\geq \left[ \int_{s_0}^{T} \left( \frac{d\hat{\mathbf{x}}_0}{dt} - \frac{d\mathbf{x}}{dt} \right) dt \right]^2$$
$$= [(\mathbf{x}(T) - \hat{\mathbf{x}}_0(T)) - (\mathbf{x}(s_0) - \hat{\mathbf{x}}_0(s_0))]^2.$$

So

$$\|\mathbf{x} - \hat{\mathbf{x}}_0\|_\infty = |\mathbf{x}(T) - \hat{\mathbf{x}}_0(T)|$$
$$\leq |\mathbf{x}(s_0) - \hat{\mathbf{x}}_0(s_0)| + \sqrt{J(\hat{\mathbf{x}}_0)}$$
$$\leq \beta_1 \kappa \left\| \frac{d^2\hat{\mathbf{x}}_0}{dt^2} \right\|_{L^2[0,T]} |\tau|^{1/2} + \beta_2 \kappa |\tau|$$
$$+ \kappa(4\sqrt{6\kappa} + \beta_3)\beta_4\sqrt{T}|\tau|^{3/2} + \beta_5\sqrt{T}|\tau|^3 + \beta_6|\tau|^{7/2},$$

where

$$\beta_5 = \sqrt{(4K^2 + 2)}C_1 K_4.$$

- If $\hat{\mathbf{x}}_0(s_0) - \mathbf{x}(s) < 0$, then there exists a $s' \in (s_0, T)$ such that $\hat{\mathbf{x}}_0(s') - \mathbf{x}(s') = 0$. So we can still use the same argument in the first case with the lower limits $s_0$ of all the integrals replaced by $s'$.  $\square$

**Acknowledgments.** The authors are very grateful to the editor, Professor Susan Murphy, the Associate Editor and two referees for their helpful comments and suggestions that greatly improved the paper. We thank Professor Jim Ramsay for helpful discussion.

Yale School of Public Health
Yale University
P.O. Box 208034, 60 College Street
New Haven, Connecticut 06520-8034
USA
E-mail: xin.qi@yale.edu
        hongyu.zhao@yale.edu